\documentclass[12pt]{article}

\usepackage{amsmath}
\usepackage{amssymb}
\usepackage[mathscr]{eucal}
\usepackage{graphicx}
\usepackage{color}

\usepackage{bm}
\usepackage{mathrsfs}
\usepackage{amsthm}
\usepackage{enumerate}
\usepackage[mathscr]{eucal}
\usepackage{eqlist}
\usepackage{graphicx}

\newtheorem{Theorem}{\sc Theorem}
\newtheorem{Definition}[Theorem]{\sc Definition}
\newtheorem{Proposition}[Theorem]{\sc Proposition}
\newtheorem{Lemma}[Theorem]{\sc Lemma}
\newtheorem{Corollary}[Theorem]{\sc Corollary}
\newtheorem{Remark}[Theorem]{\sc Remark}

\newtheorem{Problem}[Theorem]{\sc Problem}

\newcommand{\R}{{\if mm {\rm I}\mkern -3mu{\rm R}\else \leavevmode
		\hbox{I}\kern -.17em\hbox{R} \fi}}

\def\sqr#1#2{{
		\vcenter{
			\vbox{\hrule height.#2pt
				\hbox{\vrule width.#2pt height#1pt \kern#1pt
					\vrule width.#2pt
				}
				\hrule height.#2pt
			}
		}
}}

\def\real{\mathbb{R}}

\def\lista#1
{{ \itemindent 0.0cm \labelsep .2cm \leftmargin 0.8cm \rightmargin
		0.0cm \labelwidth 0.6cm \topsep 0.0mm
		\parsep 0.0mm
		\itemsep 0.0mm
		\begin{list}{}
			{ \setlength{\leftmargin}{.8cm} \setlength{\rightmargin}{0.0cm}
				\setlength{\parsep}{0.0mm} \setlength{\topsep}{.0mm}
				\setlength{\parskip}{.0cm} \setlength{\itemsep}{.0cm} }
			{#1}\end{list}} }

\begin{document}
\title{
A Class of Generalized Mixed Variational-Hemivariational Inequalities I: \\
Existence and Uniqueness Results
	\thanks{
	\, This project has received funding from the European Union's Horizon 2020 Research and Innovation Programme under the Marie Sk{\l}odowska-Curie grant agreement No. 823731 -- CONMECH. It is supported by the National Science Center of Poland under Maestro Project
No. UMO-2012/06/A/ST1/00262, and
National Science Center of Poland under Preludium Project No. 2017/25/N/ST1/00611.
The second author is also supported by the  Natural Sciences Foundation of Guangxi Grant No. 2018JJA110006, the Beibu Gulf University Project No. 2018KYQD03,
and the Project Financed by the Ministry of Science and Higher Education of Republic of Poland under Grant No.
4004/GGPJII/H2020/2018/0.
}}
	
\author{Yunru Bai~\footnote{\,Jagiellonian University in Krakow, Faculty of Mathematics and Computer Science, ul. Lojasiewicza 6, 30348 Krakow, Poland.
E-mail address:	yunrubai@163.com.},\ \ \
Stanis{\l}aw Mig\'orski~\footnote{\,College of Sciences, Beibu Gulf University, Qinzhou, Guangxi 535000, P.R. China, and Jagiellonian University in Krakow, Chair of Optimization and Control, ul. Lojasiewicza 6, 30348 Krakow, Poland.
Tel.: +48-12-6646666. E-mail address: stanislaw.migorski@uj.edu.pl.}
\ \ \mbox{and} \ \
Shengda Zeng~\footnote{\,Jagiellonian University in Krakow, Faculty of Mathematics and Computer Science, ul. Lojasiewicza 6, 30348 Krakow, Poland. Corresponding author. Tel.: +86-18059034172.
E-mail address:	shengdazeng@gmail.com; shdzeng@hotmail.com; zengshengda@163.com.
}}
	
\date{}
\maketitle
	
\noindent {\bf Abstract.} \
We investigate a generalized Lagrange multiplier system in a Banach space, called a mixed variational-hemivariational inequality
(MVHVI, for short),
which contains a hemivariational inequality
and a variational inequality.
First, we employ the Minty technique
and a monotonicity argument to establish
an equivalence theorem, which provides three different equivalent formulations of the inequality problem. Without compactness for one of operators in the problem,
a general existence theorem
for (MVHVI) is proved by using the Fan-Knaster-Kuratowski-Mazurkiewicz principle combined with methods of nonsmooth analysis. Furthermore,
we demonstrate several crucial properties of the solution set to (MVHVI) which include boundedness, convexity, weak closedness,
and continuity.
Finally, a uniqueness result with respect to the first component of the solution for the inequality problem is proved by using the Ladyzhenskaya-Babu\v{s}ka-Brezzi (LBB) condition. All results are obtained in a general functional framework in reflexive Banach spaces.

\smallskip
	
	
\smallskip
	

\section{Introduction and problem  statement}\label{Introduction}

In many complicated physical processes and engineering applications, mathematical models
based on variational inequality formulations and their generalizations play an important role.
Recently, a new class of systems with Lagrange multipliers which consists of a variational inequality has drawn a great attention.
The reason is that such systems are a powerful mathematical tool to model and solve a variety of problems in engineering areas such as dynamic vehicle routing problems, contact problems in mechanics, the behavior of Navier-Stokes fluids, the penetration phenomenon of magnetic field,
etc.

The most representative recent results in this area are the following: Cojocaru-Matei~\cite{matei2018narwa}
who have discussed the unique solvability for a class of frictional contact problems governed by the $p$-Laplace operator, which can be formulated as a mixed variational inequality;
Matei et al.~\cite{matei2018aa}
have employed the Lagrange multipliers method to consider a deformable body in frictionless unilateral contact with a moving rigid obstacle, and explored an efficient algorithm approximating the weak solution for a more general case
of a two-body contact problem including friction; Han-Reddy~\cite{hanreddy}
who have analyzed the finite element method
for a class of mixed variational inequalities of the second kind which arises in elastoplastic problems; Sofonea-Matei~\cite{sofoneamatei2015jogo} who have considered a new class of mixed variational problems, and proved existence, uniqueness as well as continuous dependence results by applying generalized saddle point formulations
and various estimates, combined with a fixed point argument.
We refer the reader  to~\cite{Matei,matei2014rwa,matei2014mma,matei2015,matei2018,matei2010,matei2018aa,matei2016ima,SM2,SM} and the references therein for a more detailed discussion of this topic.

On the other hand, the notion of
a hemivariational inequality was first introduced and studied by P.D. Panagiotopoulos~\cite{pana1983,Pa1,Pa2}
in the early 1980s who used this mathematical tool to describe and solve complicated problems modeling various physical phenomena.
After that, more and more researchers are attracted to boost the development of the theory and applications of hemivariational inequalities, since they can be applied to a wide range of engineering problems involving nonmonotone and possibly multivalued constitutive and interface laws for deformable bodies, see e.g.~\cite{han1,bartosz,carl,han2,han3,smo1,migorskiochal2009,migorskizeng1,migorskizeng2,migorskizeng3,motreanu1,SMHistory}.
Very recently, Matei~\cite{matei2018}
has studied an abstract system with Lagrange multipliers, called a mixed
variational-hemivariational inequality, which consists of a hemivariational inequality and a variational inequality, and then demonstrated three existence theorems which are
illustrated by two applications.
However, in paper~\cite{matei2018},
some problems concerning mixed variational-hemivariational inequalities,
such as uniqueness, are left open.
Based on this motivation, in this paper,
we will develop a new class of abstract mixed variational-hemivariational inequalities
in a general functional framework.

Let $V$ and $E$ be reflexive Banach spaces, and $\Lambda$ be a nonempty subset of $E$.
We denote by
$\langle \cdot, \cdot \rangle$
the duality pairing between $V$ and its dual $V^*$.
Let $X$ be another Banach space.
Given an operator $A\colon V\to V^*$,
a function $J\colon X\to \real$,
a bilinear function $b\colon V\times E\to \real$, an operator $\gamma \colon V \to X$
and an element $f\in V^*$,
the purpose of this paper is to study the following abstract generalized mixed variational-hemivariational inequality.
\begin{Problem}\label{problem1}
Find $(u,\lambda)\in V\times \Lambda$
such that the following two inequalities hold
\begin{eqnarray}
&&
\hspace{-0.5cm}
\langle A(u),v-u\rangle+b(v-u,\lambda)+J^0(\gamma u;\gamma v-\gamma u) \ge \langle f,v-u\rangle  \ \, \mbox{\rm for all}\ \ v\in V, \label{eqns1}\\[2mm]
&&\hspace{-0.5cm}
b(u,\rho-\lambda)\le 0
\ \ \ \mbox{\rm for all}\ \
\rho\in \Lambda.\label{eqns2}
\end{eqnarray}
\end{Problem}

To highlight the motivation to study  Problem~\ref{problem1},
we mention below its particular cases.

\smallskip

{\rm (i)}
Let $j\colon \real^d\to \real$ be a Lipschitz continuous function and $\Omega\subset \real^d$ ($d\ge 2$) be a bounded domain with smooth boundary $\Gamma=\partial \Omega$.
If $\gamma\colon V\to L^p(\Gamma;\real^r)$
is a linear, bounded and compact operator
with $2\le p<\infty$ and $r\ge1$,
and $J$ is defined by
\begin{equation*}
J(v)=\int_{\Gamma_1}j(\gamma v(x))\,d\Gamma
\ \ \ \mbox{\rm for all}\ \ v\in V,
\end{equation*}
where $\Gamma_1\subset \Gamma$ is such that $\mbox{meas}(\Gamma_1)>0$,
then Problem~\ref{problem1} reduces to
\begin{eqnarray*}
&&
\hspace{-0.6cm}
\displaystyle\langle A(u),v-u\rangle+b(v-u,\lambda)
+\int_{\Gamma_1}j^0(u(x);v(x)-u(x))\,d\Gamma
\ge \langle f,v-u\rangle
\ \ \mbox{\rm for all}\ \ v\in V, \\
&&\hspace{-0.6cm}
b(u,\rho-\lambda)\le 0
\ \ \ \mbox{\rm for all}\ \ \rho\in \Lambda,
\end{eqnarray*}
which has been recently studied by Matei~\cite{matei2018}.

\smallskip

{\rm (ii)}
If  $J\equiv0$, then Problem~\ref{problem1}
becomes
\begin{equation*}
\left\{\begin{array}{lll}
\langle A(u),v\rangle+b(v,\lambda)= \langle f,v\rangle
\ \ \mbox{\rm for all} \ \ v\in V, \\[1mm]
b(u,\rho-\lambda)\le 0 \ \ \mbox{\rm for all}  \ \ \rho\in \Lambda.
\end{array}\right.\\
\end{equation*}
This mixed variational inequality has been investigated by Cojocaru-Matei~\cite{matei2018narwa}.

\medskip

The aim of this paper is to extend the theoretical results from~\cite{matei2018} to a genera\-li\-zed mixed variational-hemivariational inequality in a general functional framework, Problem~\ref{problem1}, and provide positive answers to open problems remained in~\cite{matei2018}.
The main novelties of the paper are described as follows.

First, in the study of Problem~\ref{problem1}, we do not require that function $J$ is Lipschitz continuous and the operator $\gamma\colon V\to X$ is compact. This extends the scope of applications for mixed variational-hemivariational inequality.
Besides, the main core of the proof is completely different from the one carried out in~\cite{matei2018}, here we employ the well-known Fan-Knaster-Kuratowski-Mazurkiewicz theorem, not a fixed point principle.

Second, our results can be applied to a special case of Problem~\ref{problem1} in which  $b\equiv0$, and the problem reduces to the following ``pure" hemivariational inequality
\begin{equation*}
\langle A(u),v-u\rangle +J^0(\gamma u;\gamma v-\gamma u) \ge \langle f,v-u\rangle  \ \ \mbox{\rm for all}\ \ v\in V.
\end{equation*}
In fact, the above inequality has been explored by many scholars from the mathematical and application points of view
under the crucial hypothesis that operator $\gamma$ is compact, see e.g.~\cite{costearadulescuJOGO2012,tanghuang2013JOGO,xiaofanghuangJOGO2015}. However, in our results,
we will overcome this assumption.

Third, for the first time, we provide the uniqueness theorem to Problem~\ref{problem1}
with respect to the first component $u \in V$
of solution.
In the meanwhile, we develop several important properties
of the solution set to Problem~\ref{problem1},
which include boundedness, convexity, weak closedness, continuity, etc.
We believe that those results will be found useful in a number of complex problems involving a mixed variational-hemivariational inequality
as a subsystem, for instance,
in optimal control problems driven by mixed variational-hemivariational inequalities.

The outline of the paper is as follows.
Basic notation and preliminary material
needed in the sequel are recalled in Section~\ref{Preliminaries}.
In Section~\ref{Section3}, we deliver our main results concerning Problem~\ref{problem1} which include a Minty type equivalence result,
a general existence theorem,
several significant properties of the solution set,
and a uniqueness result.

\section{Background material}
\label{Preliminaries}

In this section, we briefly review basic notation and some results which are needed in the sequel. For more details, we refer to monographs~\cite{cf,DMP1,DMP2,Zeidler}.

Throughout the paper, we denote by
$\langle \cdot, \cdot \rangle_{Y^*\times Y}$
the duality pairing between a Banach
space $Y$ and its dual $Y^*$.
The norm in a normed space $Y$ is denoted by
$\| \cdot \|_Y$.
Given a subset $D$ of $Y$, we write
$\|D\|_Y = \sup \{\|v\|_Y \mid v \in Y \}$.
If no confusion arises, we often drop the subscripts.
Besides, we denote by $\mathcal{L}(Y_1, Y_2)$
the space of linear and bounded operators
from a normed space $Y_1$ to a normed space
$Y_2$ endowed with the usual norm $\|\cdot\|_{\mathcal{L}(Y_1, Y_2)}$.

We begin with definitions and properties of semicontinuous multivalued mappings.
\begin{Definition}\label{defusc}
Let $X$ and $Y$ be topological spaces, and $F\colon X\rightarrow 2^Y$ be a multivalued mapping.  We say that $F$ is

{\rm{(i)}} upper semicontinuous (u.s.c., for short) at $x\in X$ if, for
every open set $O\subset Y$

\quad \ with $F(x)\subset O$ there exists a
neighborhood $N(x)$ of $x$ such that $F(N(x)):=$

\quad \ $\cup_{y\in
N(x)}F(y)\subset O$. If this holds for every $x\in X$, then $F$
is called upper semi-

\quad \
continuous.

{\rm{(ii)}} closed at $x_0\in X$, if for every sequence $\{(x_n,y_n)\}\subset \mbox{Gr}(F)$ such that $(x_n,y_n)\to$

\quad \ $(x_0,y_0)$ in $X\times Y$, we have $(x_0,y_0)\in \mbox{Gr}(F)$, where $\mbox{Gr}(F)$ is the graph of the

\quad \ multivalued mapping $F$ defined by
    \begin{eqnarray*}
    \mbox{Gr}(F):=\{(x,y)\in X\times Y\,\mid\,y\in F(x)\}.
    \end{eqnarray*}

\quad \ We say that $F$ is closed (or $F$ has a closed graph), if it is closed at every

\quad \ $x_0\in X$.
\end{Definition}

The following theorem gives a criterium
for upper semicontinuity.
\begin{Proposition}\cite[Proposition 3.8]{smo1} \label{posusc}
Let $X$ and $Y$ be two topological spaces, and $F\colon X\rightarrow 2^Y$. The
following statements are equivalent:

{\rm{(i)}} $F$ is u.s.c..

{\rm{(ii)}} for every closed set $C\subset Y$, the set
$$
F^{-}(C):=\{x\in X \mid
F(x)\cap C\neq\emptyset\},
$$

\quad\ is closed in $X$.
\end{Proposition}

\begin{Theorem}
\cite[Proposition~4.1.9]{DMP1}\label{closedthem}
Let $X$ be a topological space, $Y$ be a regular topological space, and $F\colon X\to 2^Y$
be an upper semicontinuous multivalued mapping
with closed values. Then $F$ is closed.
\end{Theorem}

Let $(V,\|\cdot\|_V)$ be a Banach space. A function $J\colon V\to \real$ is called to be locally Lipschitz continuous at $u\in V$, if there exist a neighborhood $N(u)$ of $u$ and a constant $L_u>0$ such that
\begin{equation*}
|J(w)-J(v)|\le L_u\|w-v\|_V
\ \ \mbox{for all}\ \ w, v\in N(u).
\end{equation*}

\begin{Definition}\label{SUB}
Given a locally Lipschitz function
$J \colon V \to \real$, we denote by
$J^0 (u; v)$ the generalized (Clarke) directional derivative of $J$ at the point $u\in V$ in the direction $v\in V$ defined
by
\begin{equation*}\label{defcalark}
J^0(u;v) = \limsup
\limits_{\lambda\to 0^{+}, \, w\to u} \frac{J(w+\lambda v)-J(w)}{\lambda}.
\end{equation*}
The generalized gradient of $J \colon V \to \mathbb{R}$ at $u\in V$ is given by
\begin{equation*}
\partial J(u) = \{\, \xi\in V^{*} \mid J^0 (u; v)\ge
\langle\xi, v\rangle \ \ \mbox{\rm for all} \ \ v \in V \, \}.
\end{equation*}
\end{Definition}

The generalized gradient and generalized directional derivative of a locally Lipschitz function enjoy many nice properties and rich calculus. Here we just collect below some basic
and crucial results,
see e.g.~\cite[Proposition 3.23]{smo1}.
\begin{Proposition}\label{clarkesub}
Let $J\colon V\to \real$ be a locally Lipschitz continuous function.
Then

\smallskip

{\rm (i)}
for each $u\in V$, the function $V\ni v\mapsto J^0(u;v)\in \real$ is positively homogeneous,

\quad \
subadditive, and  satisfies
$|J^0(u;v)|\le L_u\|v\|_V$
for all $v \in V$, where $L_u>0$ is

\quad \ the Lipschitz constant of $J$ near $u$.

\smallskip

{\rm (ii)}
the function $V \times V \ni (u,v)\mapsto J^0(u;v)
\in \real$ is upper semicontinuous.

\smallskip

{\rm (iii)}
for each $v\in V$, we have
$J^0(u;v)=\max\big\{\langle u^*,v\rangle \mid u^*\in\partial J(u)\big\}$.
\end{Proposition}

We conclude this section with the following Fan-Knaster-Kuratowski-Mazurkiewicz theorem (F-KKM theorem, for short) which will play an important role in the proof of existence of solutions to the inequality problems in Section~\ref{Section3}.
Its proof can be found in Ky Fan~\cite{fk}.

\begin{Theorem}\cite{fk}\label{kkm}
Let $K$ be a nonempty subset of a Hausdorff topological vector space $E$ and
$G \colon K\to  2^E$ be a multivalued mapping with the following properties:

{\rm (a)} \
$G$ is a KKM mapping, that is, for any $\{v_{1},v_{2},\ldots,v_{N}\}\subset K$,
one has that

\qquad  its convex hull $co\{v_{1},v_{2},\ldots,v_{n}\}$
is contained in $\bigcup_{i=1}^{n}G(v_{i})$,

{\rm (b)} \
for every $v\in K$, $G(v)$ is closed in $E$,

{\rm (c)} \
for some $v_{0}\in K$,
$G(v_{0})$ is compact in $E$.

\smallskip

\noindent
Then, we have $\bigcap_{v\in K}G(v)\neq\emptyset$.
\end{Theorem}

\section{Existence and uniqueness results}\label{Section3}

The section is devoted to deliver the main results of this paper, which contain five theorems and two corollaries.
More precisely, the first theorem, Theorem~\ref{theorem1}, provides three various equivalent formulations for Problem~\ref{problem1} by using the Minty approach and a monotone argument.
In the second theorem, Theorem~\ref{theorem2},
we employ the Fan-Knaster-Kuratowski-Mazurkiewicz theorem and the theory of nonsmooth analysis to establish an existence result to Problem~\ref{problem1}, in which we do not require that the operator $\gamma$ is compact.
Next result, Theorem~\ref{theorem3} is devoted to explore some important properties of solution set of Problem~\ref{problem1}, which include boundedness, convexity, weak closedness,
and continuity.
Subsequently, a uniqueness theorem, Theorem~\ref{theorem4}, for Problem~\ref{problem1} is established by using the Ladyzhenskaya-Babu\v{s}ka-Brezzi (LBB) condition.
The last result, Theorem~\ref{theorem5}, presents a continuity result (or stability result) for the solution mapping.

\medskip

To establish main results on Problem~\ref{problem1}, we now impose the following assumptions on its data.
Let $(V,\|\cdot\|_V)$ and $(E,\|\cdot\|_E)$ be two reflexive Banach spaces.

\medskip
\noindent
${\underline{H(h)}}$: $h\colon V\to \real$ is such that $h(0_V)=0$ and

\smallskip
\noindent
(i) $\limsup_{t\to0^+}\frac{h(tv)}{t}\ge 0$
for all $v\in V$.

\smallskip
\noindent
(ii) for all $\{v_n\}\subset V$ with $v_n\to v$ weakly in $V$, we have
$h(v) \le \limsup_{n\to\infty}h(v_n)$.

\smallskip

\noindent
(iii) for all $v\in V\backslash\{0_V\}$, we have $h(v)> 0$.

\medskip
\noindent
${\underline{H(J)}}$: $J\colon X\to \real$ is such that

\smallskip

\noindent
(i) $J$ is locally Lipschitz continuous.

\smallskip

\noindent
(ii) there exist $\theta\ge0$, $\alpha_J\ge0$
and $\beta_J>0$ such that
\begin{equation*}
J^0(v; -v)\le \alpha_J+\beta_J\|v\|_X^\theta
\ \ \ \mbox{\rm for all}\ \ v\in X.
\end{equation*}

\smallskip

\noindent
(iii) the multivalued mapping
$X\ni v \mapsto\partial J(v)\subset X^*$
is bounded, i.e., $\partial J$ maps bounded
subsets of $X$ into bounded subsets of $X^*$.

\medskip
\noindent
${\underline{H(A)}}$: $A\colon V\to V^*$
is such that

\smallskip

\noindent
(i) for any $w$, $v\in V$ fixed, it holds
\begin{equation*}
\limsup_{t\to 0^+}\, \langle A(tw+(1-t)v),
w-v\rangle\le \langle A(v), w-v\rangle.
\end{equation*}

\smallskip
\noindent
(ii) the mapping $A\cdot + \gamma^*\partial J(\gamma\cdot)\colon V\to 2^{V^*}$ is $h$-relaxed monotone on $V$, i.e.,
\begin{equation*}
\langle A(u)+\gamma^*\xi_u-A(v)-\gamma^*\xi_v, u-v\rangle\ge h(u-v)
\end{equation*}

\ \  for all $\xi_u\in\partial J(\gamma u)$, $\xi_v\in\partial J(\gamma v)$
and $u$, $v\in V$.

\smallskip

\noindent
(iii) $A$ is coercive in the following sense
\begin{equation*}
\lim_{v\in V,\, \|v \|_V\to+\infty}
\, \frac{\langle Av, v\rangle}
{\|v\|_V^{\max\{\theta,1\}}}=+\infty,
\end{equation*}

\ \ where $\theta\geq0$ is given in hypothesis $H(J)$(ii).

\smallskip
\noindent
(iv) $A$ is a bounded operator.

\medskip
\noindent
${\underline{H(b)}}$: The bilinear function $b\colon V\times E\to \real$ is bounded and satisfies the following

\quad\  inequality
\begin{equation}\label{inequalityb}
\inf_{\rho\in E\setminus\{0_E\}}
\sup_{v\in V\setminus\{0_V\}}
\frac{b \big(v,\rho\big)}
{\|v\|_V\|\rho\|_E}\ge \alpha_b
\end{equation}

\quad\  for some $\alpha_b>0$.

\medskip

\noindent
${\underline{H(\gamma)}}$: $\gamma\colon V\to X$ is a linear and continuous operator.

\medskip

In the following we comment on the above
hypotheses.
\begin{Remark}\label{remark1}
Various kinds of monotonicity of operator
$A \cdot + \gamma^*\partial J(\gamma\cdot)$
can be obtained by choosing a suitable function $h$.	
%
%
It is worth to mention that if hypotheses $H(h)$(i) and (iii) are specified by $h(tv)=t^\tau h(v)$ for all $v\in V$ and $t>0$ with $\tau>1$, and $h(v)\ge c_h\, \|v\|_V^q$
for all $v\in V$ with some $c_h>0$ and $q>1$, respectively, then hypotheses $H(h)$ reduces
to the one considered by  Cojocaru-Matei~\cite{matei2018narwa}.
In particular, function
$h(v)=m_A \|v\|^\tau$ for $v \in V$ with
$m_A>0$ and $\tau\ge 1$ enjoys $H(h)$,
and then $H(A)(ii)$ means that $A \cdot +\gamma^*\partial J(\gamma\cdot)$ is $\tau$-strongly monotone.
Note also that hypothesis $H(A)$(i) is weaker then the hemicontinuity of operator $A$,
see~\cite[Definition 3.68]{smo1}.

Let us turn to the hypotheses $H(J)$. If the generalized gradient $\partial J$ has a sublinear growth, namely,
\begin{equation*}
\|\partial J(v)\|_{X^*} \le c_J + d_J\|v\|_X
\ \ \mbox{\rm for all}\ \ v\in X
\end{equation*}
for some $c_J\ge 0$ and $d_J>0$, then
hypothesis $H(J)$(iii) is clearly satisfied.

The inequality (\ref{inequalityb}) is usually called the Ladyzhenskaya-Babu\v{s}ka-Brezzi (LBB) condition which widely appears in the literature.
\end{Remark}

The first result of the paper provides three different equivalent formulations of Problem~\ref{problem1} by applying the Minty approach, in which $V$ and $\Lambda$ are replaced by the nonempty, closed and convex subsets of $V$ and $E$, respectively.
\begin{Theorem}\label{theorem1}
Let $K$  and $Y$ be nonempty, closed and convex subsets of $V$ and $E$, respectively. Assume that hypotheses $H(A)$(i)--(ii), $H(J)$(i), $H(h)$(i), and $H(\gamma)$ hold. If $b\colon V\times E\to \real$ is a bilinear and bounded function, then $(u,\lambda)\in K\times Y$ is a solution to the following mixed variational-hemivariational  inequality
\begin{eqnarray}
&&
\langle A(u),v-u\rangle+b(v-u,\lambda)+J^0(\gamma u;\gamma v-\gamma u) \ge \langle f,v-u\rangle  \ \, \mbox{\rm for all}\ \ v\in K, \label{eqn1}\\[2mm]
&&
b(u,\rho-\lambda)\le 0
\ \ \ \mbox{\rm for all}\ \
\rho\in Y,\label{eqn2}
\end{eqnarray}
if and only if it solves one of the following problems

\smallskip

{\rm (i)}
$(u,\lambda)\in K\times Y$ is such that
\begin{eqnarray}
&&
\hspace{-1.0cm}
\left\{\begin{array}{lll}\label{eqns3}
\langle Av, v-u\rangle
+b(v-u,\lambda)+J^0(\gamma v;\gamma v-\gamma u) \\[1mm]
\qquad\quad\ge \langle f,v-u\rangle+h(v-u) \ \ \mbox{\rm for all}\ \ v\in K,
\end{array}\right.\\[2mm]
&&
\hspace{-0.3cm}
b(u,\rho-\lambda)\le 0\ \ \ \mbox{\rm for all}\ \
\rho\in Y. \label{eqns4}
\end{eqnarray}

\smallskip

{\rm (ii)}
$(u,\lambda)\in K\times Y$ is such that
\begin{eqnarray}
&&
\hspace{-1.0cm}
\left\{\begin{array}{lll}\label{eqns5}
\langle Au,v-u\rangle+b(v,\lambda)-b(u,\rho) +J^0(\gamma u;\gamma v-\gamma u)\\[2mm]
\qquad\quad
\ge \langle f,v-u\rangle
\ \ \mbox{\rm for all}\
v\in K\ \mbox{\rm and all}\ \rho\in Y.
\end{array}\right.
\end{eqnarray}

\smallskip
{\rm (iii)}
$(u,\lambda)\in K\times Y$
is such that
\begin{eqnarray}
&&
\hspace{-1.0cm}
\left\{\begin{array}{lll}\label{eqns6}
\langle Av,v-u\rangle
+b(v,\lambda)-b(u,\rho)+J^0(\gamma v;\gamma v-\gamma u)\\[2mm]
\qquad\quad \ge \langle f,v-u\rangle+h(v-u)
\ \ \mbox{\rm for all}\
v\in K\ \mbox{\rm and all}\ \rho\in Y.
\end{array}\right.
\end{eqnarray}
\end{Theorem}

\noindent
{\bf Proof.}
(i) Let $(u,\lambda)\in K\times Y$ be a solution to problem (\ref{eqn1}) and (\ref{eqn2}).
It is obvious that inequality (\ref{eqn2}) coincides with (\ref{eqns4}).
Moreover, the $h$-relaxed monotonicity of operator
$A \cdot +\gamma^*\partial J(\gamma\cdot)$
leads to
\begin{equation*}
\langle Av,v-u\rangle+\langle\xi_v,\gamma(v-u)\rangle_{X^*\times X} \ge \langle Au,v-u\rangle+\langle\xi_u,\gamma(v-u)\rangle_{X^*\times X}
+h(v-u)
\end{equation*}
for all $\xi_u\in\partial J(\gamma u)$, $\xi_v\in\partial J(\gamma v)$ and all
$u$, $v \in V$.
Taking into account the above inequality,
the property
\begin{equation*}
J^0(\gamma u;\gamma v-\gamma u)=\max_{\xi_u\in \partial J(\gamma u)}\langle\xi_u,\gamma(v-u)\rangle_{X^*\times X}=\langle\widetilde\xi_u,\gamma(v-u)\rangle_{X^*\times X}
\end{equation*}
for some
$\widetilde\xi_u\in \partial J(\gamma u)$,
and inequality (\ref{eqn1}), we obtain
\begin{eqnarray*}
&&\hspace{-0.5cm}
h(v-u)+\langle f,v-u\rangle \\[2mm]
&&\le \langle Au,v-u\rangle+J^0(\gamma u;\gamma v-\gamma u)+b(v-u,\lambda)+h(v-u)\\[2mm]
&&= \langle Au,v-u\rangle+\langle\widetilde\xi_u,\gamma(v-u)\rangle_{X^*\times X}+b(v-u,\lambda)+h(v-u)\\[2mm]
&&\le \langle Av,v-u\rangle+\langle\xi_v,\gamma(v-u)\rangle_{X^*\times X}+b(v-u,\lambda)\ \,\mbox{ (for all $\xi_v\in\partial J(\gamma v)$)}\\[2mm]
&&\le \langle Av,v-u\rangle+J^0(\gamma v;\gamma v-\gamma u)+b(v-u,\lambda)
\end{eqnarray*}
for all $v\in K$. So, $(u,\lambda)\in K\times Y$
is also a solution to problem (\ref{eqns3}) and (\ref{eqns4}).

Conversely, let $(u,\lambda)\in K\times Y$ be a solution to problem (\ref{eqns3}) and (\ref{eqns4}).
Then, (\ref{eqn2}) holds due to (\ref{eqns4}).
It is enough to obtain (\ref{eqn1}).
Let $w\in K$,  and $t\in(0,1)$ be arbitrary.
Taking $v=v_t:= tw+(1-t)u$  in (\ref{eqns3}), we employ hypotheses $H(b)$, $H(J)$(i), $H(\gamma)$, and Proposition~\ref{clarkesub} to get
\begin{eqnarray*}
&&
\hspace{-0.5cm}
h(t(w-u))+t\langle f,w-u\rangle\\[2mm]
&&\le t\langle Av_t,w-u\rangle+tb(w-u,\lambda)+J^0(\gamma u;t(\gamma w-\gamma u))\\[2mm]
&&=t\big[\langle Av_t,w-u\rangle+b(w-u,\lambda)+J^0(\gamma u;\gamma w-\gamma u)\big],
\end{eqnarray*}
and hence
\begin{eqnarray*}
&&
\hspace{-0.5cm}
\langle Av_t,w-u\rangle+b(w-u,\lambda)
+J^0(\gamma u;\gamma w-\gamma u)\\[2mm]
&&\ge\langle f,w-u\rangle+\frac{h(t(w-u))}{t}.
\end{eqnarray*}
Passing to the upper limit as $t\to 0^+$, we now apply conditions $H(A)$(i) and $H(h)(i)$
to obtain inequality (\ref{eqn1}). This means that $(u,\lambda)\in K\times Y$ also solves problem (\ref{eqn1}) and (\ref{eqn2}).

\smallskip

(ii) Assume that $(u,\lambda)\in K\times Y$ is a solution to problem (\ref{eqn1}) and (\ref{eqn2}). The inequality (\ref{eqns5}) can be obtained readily by multiplying inequality (\ref{eqn2}) by $-1$ and then summing it up with (\ref{eqn1}).

Conversely, if $(u,\lambda)\in K\times Y$ is a solution to problem (\ref{eqns5}), then inequalities (\ref{eqn1}) and (\ref{eqn2}) are a direct consequence of (\ref{eqns5}) via inserting $\rho=\lambda$ and $v=u$ into (\ref{eqns5}), respectively.

\smallskip

(iii) From assertion (i),
it remains to show that $(u,\lambda)\in K\times Y$ is a solution to problem (\ref{eqns3}) and (\ref{eqns4}) if and only if it solves problem (\ref{eqns6}).
Indeed,
if $(u,\lambda)\in K\times Y$ is a solution to (\ref{eqns3}) and (\ref{eqns4}), then the inequality (\ref{eqns6}) is obtained easily via multiplying inequality (\ref{eqns4}) by $-1$ and by summing up the resulting inequality with (\ref{eqns3}).

For the converse, we put $\rho=\lambda$ and
$v=u$ in (\ref{eqns6}), respectively, to obtain  inequalities (\ref{eqns3}) and (\ref{eqns4}).
\hfill$\Box$

\medskip

The following theorem delivers a crucial existence result for Problem~\ref{problem1} without any compactness hypothesis on operator $\gamma$.
\begin{Theorem}\label{theorem2}
Let $\Lambda$ be a nonempty, closed and convex subset of $E$ with $0_E\in \Lambda$. If hypotheses $H(A)$, $H(b)$, $H(J)$, $H(h)$(i)--(ii), and $H(\gamma)$ are satisfied, then Problem~\ref{problem1} has at least one solution $(u,\lambda)\in V\times \Lambda$.
\end{Theorem}

\noindent
{\bf Proof.}
It follows from Theorem~\ref{theorem1} that it is enough to prove that problem (\ref{eqns5}) admits a solution with $K=V$ and $Y=\Lambda$.
The proof will be carried out in three steps.

Let $r$, $s\in \mathbb N$.
We define sets $K(r)$ and $Y(s)$ by
\begin{equation*}
K(r):= \overline B_{V}(r)\subset V\quad\mbox{ and }\quad Y(s):=\overline B_{E}(s)\cap \Lambda\subset \Lambda,
\end{equation*}
where $\overline B_{H}(r)$ stands for a closed ball with centre $0_H$ and radius $r>0$ in a space $H$.

\medskip

\noindent
{\bf Step 1.}
{\it For all $r$, $s\in\mathbb N$ fixed, the following problem admits a solution $(u_r,\lambda_s)\in K(r)\times Y(s)$ such that}
\begin{eqnarray}
&&
\hspace{-1.0cm}
\left\{\begin{array}{lll}\label{eqnss7}
\langle Au_r,v-u_r\rangle+b(v,\lambda_s)-b(u_r,\rho) +J^0(\gamma u_r;\gamma v-\gamma u_r)\\[2mm]
\qquad\quad
\ge \langle f,v-u_r\rangle
\ \ \mbox{\rm for all}\
v\in K(r)\ \mbox{\rm and all}\ \rho\in Y(s).
\end{array}\right.
\end{eqnarray}

\smallskip

We consider a multivalued mapping
$G\colon K(r)\times Y(s)\to 2^{K(r)\times Y(s)}$ given by
\begin{eqnarray}\label{eqnss8}
&&G(v,\rho):=\big\{(u,\lambda)\in K(r)\times Y(s)\,\mid\,\langle Av,v-u\rangle+b(v,\lambda)-b(u,\rho)\nonumber\\[2mm]
&&\hspace{2cm}\quad+\inf_{\xi_v\in\partial J(\gamma v)}\langle\xi_v,\gamma(v-u)\rangle_{X^*\times X}\ge \langle f,v-u\rangle+h(v-u)\big\}.
\end{eqnarray}
It is obvious that for each $(v,\rho)\in K(r)\times Y(s)$, the set $G(v,\rho)$ is nonempty, since $(v,\rho)\in G(v,\rho)$.
We now demonstrate that for each
$(v,\rho)\in K(r)\times Y(s)$ fixed,
the set $G(v,\rho)$ is weakly closed. Let $\{(u_n,\lambda_n)\}\subset G(v,\rho)$ be a sequence such that $u_n\to u$ weakly in $V$ and $\lambda_n\to \lambda $ weakly in $E$.
We have
\begin{eqnarray*}
&&
\langle Av,v-u_n\rangle+b(v,\lambda_n)-b(u_n,\rho)
+\langle\xi_v,\gamma(v-u_n)
\rangle_{X^*\times X}\\[2mm]
&&
\qquad\ge \langle f,v-u_n\rangle+h(v-u_n)
\end{eqnarray*}
for all $\xi_v\in\partial J(\gamma v)$.
Passing to the upper limit, as $n\to \infty$,
in the above inequality, from $H(h)$(ii),
one has
\begin{eqnarray*}
&&
\hspace{-0.5cm}
\langle f,v-u\rangle+h(v-u)\\[2mm]
&&\le \lim_{n\to \infty}\langle f,v-u_n\rangle+\limsup_{n\to \infty}h(v-u_n)\\[2mm]
&&= \limsup_{n\to \infty}\big[\langle f,v-u_n\rangle+h(v-u_n)\big]\\[2mm]
&&\le \limsup_{n\to \infty}\bigg[\langle Av,v-u_n\rangle+b(v,\lambda_n)-b(u_n,\rho)+\langle\xi_v,\gamma(v-u_n)\rangle_{X^*\times X}\bigg]\\[2mm]
&&= \langle Av,v-u\rangle+b(v,\lambda)-b(u,\rho)+\langle\xi_v,\gamma(v-u)\rangle_{X^*\times X}
\end{eqnarray*}
for all
$\xi_v\in\partial J(\gamma v)$.
Hence,
\begin{eqnarray*}
&&
\hspace{-0.5cm}
\langle f,v-u\rangle+h(v-u)\\[2mm]
&&\le \langle Av,v-u\rangle+b(v,\lambda)-b(u,\rho)+\inf_{\xi_v\in\partial J(\gamma v)}\langle\xi_v,\gamma(v-u)\rangle_{X^*\times X}.
\end{eqnarray*}
This means that $(u,\lambda)\in G(v,\rho)$, therefore, the set $G(v,\rho)$ is weakly closed. Note that
$K(r)\times Y(s)$ is a nonempty, bounded, closed, and convex subset of $V\times E$.
Moreover, from the reflexivity of $V\times E$, it follows that the set $K(r)\times Y(s)$
is weakly compact and convex.
This ensures that $G(v,\rho)$ is relatively weakly compact in $V\times E$
for all $(v,\rho)\in V\times E$.
Thus, we conclude that for all
$(v,\rho)\in V\times E$,
the set $G(v,\rho)$ is weakly compact, owning to the weak closedness of $G(v,\rho)$.

Now, we can distinguish two cases:
(a) $G$ is a KKM mapping, and (b) $G$ is not a KKM mapping.

If case (a) occurs, then via invoking the F-KKM principle, Theorem~\ref{kkm}, we are able to find $(u_r,\lambda_s)\in K(r)\times Y(s)$ such that
\begin{equation*}
(u_r,\lambda_s)\in
\bigcap_{(v,\rho)\in K(r)\times Y(s)} G(v,\rho)\neq\emptyset,
\end{equation*}
that is,
\begin{equation*}
\langle Av,v-u_r\rangle+b(v,\lambda_s)-b(u_r,\rho)+ \langle\xi_v,\gamma(v-u_r)\rangle_{X^*\times X}\ge \langle f,v-u_r\rangle+h(v-u_r)
\end{equation*}
for all $\xi_v\in\partial J(\gamma v)$
and $(v,\rho)\in K(r)\times Y(s)$.
Hence, we have
\begin{equation*}
\langle Av,v-u_r\rangle+b(v,\lambda_s)-b(u_r,\rho)+ J^0(\gamma v;\gamma v-\gamma u_r)\ge
\langle f,v-u_r\rangle+h(v-u_r)
\end{equation*}
for all $(v,\rho)\in K(r)\times Y(s)$.
We now employ Theorem~\ref{theorem1} to show
that $(u_r,\lambda_s)\in K(r)\times Y(s)$
solves problem (\ref{eqnss7}).

On the other hand, when (b) holds,
then there exist $$
\{(v_1,\rho_1),(v_2,\rho_2),\ldots,(v_N,\rho_N)\}\subset K(r)\times Y(s)
$$
and $(u_r,\lambda_s)\in K(r)\times Y(s)$ with
$u_r=\sum_{i=1}^N t_iv_i$, $\lambda_s=\sum_{i=1}^Nt_i\rho_i$,
$t_i\in[0,1]$ for $i=1,2,\ldots,N$,
and $\sum_{i=1}^Nt_i=1$ such that
\begin{equation*}
(u_r,\lambda_s)\notin  \bigcup_{i=1}^NG(v_i,\rho_i).
\end{equation*}
This means that
\begin{eqnarray}\label{eqnss9}
&&\langle Av_i,v_i-u_r\rangle+b(v_i,\lambda_s)-b(u_r,\rho_i)+ \inf_{\xi_i\in\partial J(\gamma v_i)}\langle\xi_i,\gamma( v_i-u_r)\rangle_{X^*\times X}\nonumber\\[1mm]
&&\qquad< \langle f,v_i-u_r\rangle+h(v_i-u_r)
\end{eqnarray}
for $i=1,2,\ldots,N$.

\smallskip

\noindent
{\bf Claim 1.} {\it There exists a neighborhood $O$ of $(u_r,\lambda_s)$ in $V\times \Lambda$ such that whenever $(v,\rho)\in O\cap \big(K(r)\times Y(s)\big)$, there holds}
\begin{eqnarray}\label{eqnss10}
&&\langle Av_i,v_i-v\rangle+b(v_i,\rho)-b(v,\rho_i)+ \inf_{\xi_i\in\partial J(\gamma v_i)}\langle\xi_i,\gamma( v_i-v)\rangle_{X^*\times X}\nonumber\\[1mm]
&&\qquad< \langle f,v_i-v\rangle+h(v_i-v).
\end{eqnarray}

\smallskip

Arguing by contradiction, we may assume that there are $\{u_n\}\subset K(r)$, $\{\lambda_n\}\subset Y(s)$, and $\{j_n\}\subset\{1,2,\ldots,N\}$ such that $u_n\to u_r$ in $V$, $\lambda_n\to \lambda_s$ in $E$ and, for every $n\in\mathbb N$, we have
\begin{eqnarray*}
&&\langle Av_{j_n},v_{j_n}-u_n\rangle+b(v_{j_n},\lambda_n)-b(u_n,\rho_{j_n})+ \langle\xi_{j_n},\gamma( v_{j_n}-u_n)\rangle_{X^*\times X}\\[2mm]
&&\qquad\ge \langle f,v_{j_n}-u_n\rangle+h(v_{j_n}-u_n)
\end{eqnarray*}
for all $\xi_{j_n}\in\partial J(\gamma v_{j_n})$. Since
for all $n\in\mathbb N$, $j_n\in\{1,2,\ldots,N\}$,
so, without any loss of generality, we may suppose that there exists
$j_0\in\{1,2,\ldots, N\}$
such that for all $n\in\mathbb N$,
the following inequality holds
\begin{eqnarray*}
&&\langle Av_{j_0},v_{j_0}-u_n\rangle+b(v_{j_0},\lambda_n)-b(u_n,\rho_{j_0})+ \langle\xi_{j_0},\gamma( v_{j_0}-u_n)\rangle_{X^*\times X}\nonumber\\[1mm]
&&
\qquad\ge \langle f,v_{j_0}-u_n\rangle+h(v_{j_0}-u_n)
\end{eqnarray*}
for all $\xi_{j_0}\in\partial J(\gamma v_{j_0})$. If, we now pass to the upper limit, as $n\to\infty$, in the above inequality, we get
\begin{eqnarray*}
&&
\hspace{-0.5cm}
\langle f,v_{j_0}-u_r\rangle+h(v_{j_0}-u_r)\\[2mm]
&&\le \lim_{n\to\infty}\langle f,v_{j_0}-u_n\rangle+\limsup_{n\to\infty}h(v_{j_0}-u_n)\\[2mm]
&&=\limsup_{n\to\infty}\big[\langle f,v_{j_0}-u_n\rangle+h(v_{j_0}-u_n)\big]\\[2mm]
&&\le\limsup_{n\to\infty}\big[\langle Av_{j_0},v_{j_0}-u_n\rangle+b(v_{j_0},\lambda_n)-b(u_n,\rho_{j_0})+ \langle\xi_{j_0},\gamma( v_{j_0}-u_n)\rangle_{X^*\times X}\big]\\[2mm]
&&=\langle Av_{j_0},v_{j_0}-u_r\rangle+b(v_{j_0},\lambda_s)-b(u_r,\rho_{j_0})+ \langle\xi_{j_0},\gamma( v_{j_0}-u_r)\rangle_{X^*\times X}
\end{eqnarray*}
for all
$\xi_{j_0}\in \partial J(\gamma v_{j_0})$,
that is,
\begin{eqnarray*}
&&
\hspace{-0.5cm}
\langle Av_{j_0},v_{j_0}-u_r\rangle+b(v_{j_0},\lambda_s)-b(u_r,\rho_{j_0})+ \inf_{\xi_{j_0}\in \partial J(\gamma v_{j_0})}\langle\xi_{j_0},\gamma( v_{j_0}-u_r)\rangle_{X^*\times X}\\[1mm]
&&\ge\langle f,v_{j_0}-u_r\rangle+h(v_{j_0}-u_r).
\end{eqnarray*}
This is a contradiction with (\ref{eqnss9}), so, Claim 1 is valid.

\smallskip

Subsequently, from Claim 1, for every $i\in\{1,2,\ldots,N\}$,
we are able to find
$\overline\xi_i\in\partial J(\gamma v_i)$ such that
\begin{eqnarray*}
&&\langle Av_i,v_i-v\rangle+b(v_i,\rho)-b(v,\rho_i)+ \langle\overline \xi_i,\gamma( v_i-v)\rangle_{X^*\times X}\nonumber\\[1mm]
&&\qquad\le \langle f,v_i-v\rangle+h(v_i-v)
\end{eqnarray*}
for all
$(v,\rho)\in O\cap \big(K(r)\times Y(s)\big)$.
It follows from the $h$-relaxed monotonicity of operator
$A \cdot +\gamma^*\partial J(\gamma \cdot)$
that
\begin{eqnarray*}
&&
\hspace{-0.5cm}
\langle Av,v_i-v\rangle+b(v_i,\rho)-b(v,\rho_i)+ \langle \xi_v,\gamma( v_i-v)\rangle_{X^*\times X}+h(v_i-v)\nonumber\\[2mm]
&&\le \langle Av_i,v_i-v\rangle+b(v_i,\rho)-b(v,\rho_i)+ \langle\overline \xi_i,\gamma( v_i-v)\rangle_{X^*\times X}\nonumber\\[2mm]
&&\le \langle f,v_i-v\rangle+h(v_i-v)
\end{eqnarray*}
for all $\xi_v\in\partial J(\gamma v)$
and  $(v,\rho)\in O\cap \big(K(r)\times Y(s)\big)$, therefore,
\begin{equation*}
\langle Av,v-v_i\rangle+b(v,\rho_i)-b(v_i,\rho)+ \langle \xi_v,\gamma( v-v_i)\rangle_{X^*\times X}\ge \langle f,v-v_i\rangle
\end{equation*}
for all $\xi_v\in\partial J(\gamma v)$ and  $(v,\rho)\in O\cap \big(K(r)\times Y(s)\big)$. Next, multiplying the above inequality by
$t_i\ge 0$, and summing up those inequalities from $i=1$ to $N$, one obtains
\begin{equation}\label{eqnss11}
\langle Av,v-u_r\rangle+b(v,\lambda_s)-b(u_r,\rho)+ \langle \xi_v,\gamma( v-u_r)\rangle_{X^*\times X}\ge \langle f,v-u_r\rangle
\end{equation}
for all $\xi_v\in\partial J(\gamma v)$ and  $(v,\rho)\in O\cap \big(K(r)\times Y(s)\big)$, where we have used the facts
$u_r=\sum_{i=1}^N t_iv_i$, $\lambda_s=\sum_{i=1}^Nt_i\rho_i$,
$t_i\in[0,1]$ for $i=1,2,\ldots,N$,
and $\sum_{i=1}^Nt_i=1$.
Assume now that $(w,\eta)\in K(r)\times Y(s)$ is arbitrary, and consider the sequence $\{(v_n,\rho_n)\}\subset V\times E$ defined by
\begin{equation*}
v_n:=\frac{1}{n}w+(1-\frac{1}{n})u_r\quad\mbox{ and }\quad \rho_n:=\frac{1}{n}\eta+(1-\frac{1}{n})\lambda_s.
\end{equation*}
It is not difficult to find $N_1\in\mathbb N$ large enough such that $(v_n,\rho_n)\in O\cap \big(K(r)\times Y(s)\big)$
for all $n\ge N_1$.
Inserting $v=v_n$ and $\rho=\rho_n$ into (\ref{eqnss11}), it reads
\begin{eqnarray*}
&&
\hspace{-0.5cm}
\frac{1}{n}\langle f,w-u_r\rangle\\[2mm]
&&\le \frac{1}{n}\langle Av_n,w-u_r\rangle+\frac{1}{n}\, b(w,\lambda_s)+(1-\frac{1}{n})b(u_r,\lambda_s)
-\frac{1}{n}b(u_r,\eta)\\[2mm]
&&\quad-(1-\frac{1}{n})\, b(u_r,\lambda_s)
+ \frac{1}{n}\langle \xi_{v_n},\gamma( w-u_r)\rangle_{X^*\times X}\ \mbox{ (for all $\xi_{v_n}\in\partial J(\gamma v_n)$)}\\[2mm]
&&\le \frac{1}{n}\big[\langle Av_n,w-u_r\rangle+b(w,\lambda_s)-b(u_r,\eta)+ J^0(\gamma v_n;\gamma(w-u_r))\big]
\end{eqnarray*}
for all $n\ge N_1$.
If we divide both sides of the above inequality
by $\frac{1}{n}$, and then pass to the upper limit, as $n\to \infty$, we get
\begin{eqnarray*}
&&
\hspace{-0.5cm}
\langle f,w-u_r\rangle\\[2mm]
&&\le \limsup_{n\to\infty}\big[\langle Av_n,w-u_r\rangle+b(w,\lambda_s)-b(u_r,\eta)+ J^0(\gamma v_n;\gamma(w-u_r))\big]\\[2mm]
&&\le \limsup_{n\to\infty}\langle Av_n,w-u_r\rangle+b(w,\lambda_s)-b(u_r,\eta)+ \limsup_{n\to\infty} J^0(\gamma v_n;\gamma(w-u_r))\\[2mm]
&&\le \langle Au_r,w-u_r\rangle+b(w,\lambda_s)-b(u_r,\eta)
+J^0(\gamma u_r;\gamma(w-u_r)).
\end{eqnarray*}
Since $w\in K(r)$ and $\eta\in Y(s)$ are  arbitrary, we conclude that $(u_r,\lambda_s)\in K(r)\times Y(s)$ is also a solution to problem (\ref{eqnss7}).

\medskip

\noindent
{\bf Step 2.} {\it For every $s\in\mathbb N$ fixed, the following problem has at least one  solution $(u_s,\lambda_s)\in V\times Y(s)$ such that}
\begin{eqnarray}
&&
\hspace{-1.0cm}
\left\{\begin{array}{lll}\label{eqnss12}
\langle Au_s,v-u_s\rangle+b(v,\lambda_s)-b(u_s,\rho) +J^0(\gamma u_s;\gamma v-\gamma u_s)\\[2mm]
\qquad\quad
\ge \langle f,v-u_s\rangle
\ \ \mbox{\rm for all}\
v\in V\ \mbox{\rm and all}\ \rho\in Y(s).
\end{array}\right.
\end{eqnarray}

\medskip

It follows from Step 1 that for any
$r$, $s\in\mathbb N$, problem (\ref{eqnss7}) admits a solution $(u_r,\lambda_s)\in K(r)\times Y(s)$.

\smallskip

\noindent
{\bf Claim 2.} {\it There exist $r_0\in\mathbb N$ and a solution $(u_{r_0},\lambda_s)\in  K(r_0)\times Y(s)$ to problem (\ref{eqnss7}) for
$r = r_0$ such that}
\begin{equation}\label{eqnss13}
\|u_{r_0}\|_V<r_0.
\end{equation}

\smallskip

Suppose that this claim is not true, so for any $r\in\mathbb N$, for each solution $(u_r, \lambda_s)\in K(r)\times Y(s)$ of problem (\ref{eqnss7}), it holds
\begin{equation*}
\|u_r\|_V=r.
\end{equation*}
Since
$0_E\in Y(s)$ and $0_V\in K(r)$ for all
$r\in \mathbb N$, we now take $v=0_V$ and $\rho=0_E$ in (\ref{eqnss7}) to obtain
\begin{eqnarray*}
&&\langle Au_r,u_r\rangle\le b(0_V,\lambda)-b(u_r,0_E) +J^0(\gamma u_r;-\gamma u_r)+\langle f,u_r\rangle\\[2mm]
&&\quad=J^0(\gamma u_r;-\gamma u_r)
+\langle f,u_r\rangle
\le \alpha_J
+\beta_J\|\gamma\|_{\mathcal L(V,X)}^\theta
\|u_r\|_V^\theta
+\|f\|_{V^*}\|u_r\|_V,
\end{eqnarray*}
and hence,
\begin{eqnarray*}
&&\frac{\langle Au_r,u_r\rangle}{\|u_r\|_V^{\max\{\theta,1\}}}\le \frac{\alpha_J}{\|u_r\|_V^{\max\{\theta,1\}}}+\frac{\beta_J\|\gamma\|_{\mathcal L(V,X)}^\theta}{\|u_r\|_V^{\max\{0,1-\theta\}}}+\frac{\|f\|_{V^*}}{\|u_r\|_V^{\max\{0,\theta-1\}}}.
\end{eqnarray*}
Passing to the limit, as $r\to\infty$, and invoking the coercivity condition $H(A)$(iii),
we get a contradiction.
This ensures that Claim 2 is true.

\smallskip

Assume now that
$(u_{r_0},\lambda_s)\in K(r_0)\times Y(s)$
is a solution to problem (\ref{eqnss7}) for $r=r_0$ such that inequality (\ref{eqnss13}) holds.
We affirm that $(u_{r_0},\lambda_s)\in V\times Y(s)$ is also a solution to problem (\ref{eqnss12}).
Let $w\in V$ and $\eta\in Y(s)$ be arbitrary,
and $t\in(0,1)$ be small enough such that $v_t=tw+(1-t)u_{r_0}\in K(r_0)$
(thanks to inequality (\ref{eqnss13})).
Putting $v=v_t$ and $\rho=\rho_t:=t\eta+(1-t)\lambda_s$ into (\ref{eqnss7}), it holds
\begin{eqnarray*}
&&
\hspace{-0.5cm}
t\langle f,w-u_{r_0}\rangle\\[2mm]
&&\le t\langle Au_{r_0},w-u_{r_0}\rangle
+t b(w,\lambda_s)+(1-t)b(u_{r_0},\lambda_s)-tb(u_{r_0},\eta)\\[2mm]
&&\quad-(1-t) b(u_{r_0},\lambda_s)
+tJ^0(\gamma u_{r_0};\gamma w-\gamma u_{r_0}) \\[2mm]
&&=t\big[\langle Au_{r_0},w-u_{r_0}\rangle+b(w,\lambda_s)-b(u_{r_0},\eta)+J^0(\gamma u_{r_0};\gamma w-\gamma u_{r_0}) \big],
\end{eqnarray*}
and
\begin{equation*}
\langle Au_{r_0},w-u_{r_0}\rangle+b(w,\lambda_s)-b(u_{r_0},\eta) +J^0(\gamma u_{r_0};\gamma w-\gamma u_{r_0}) \ge \langle f,w-u_{r_0}\rangle
\end{equation*}
for all $w\in V$ and $\eta\in Y(s)$. Consequently, $(u_{r_0},\lambda_s)\in V\times Y(s)$ is a solution to problem (\ref{eqnss12}).

\medskip

\noindent
{\bf Step 3.} {\it Problem~\ref{problem1} has at least one solution.}

\medskip

Indeed, Step 2 guarantees that for each $s\in \mathbb N$, problem (\ref{eqnss12}) admits a solution $(u_s,\lambda_s)\in V\times Y(s)$.

\smallskip

\noindent
{\bf Claim 3.} {\it The sequences $\{u_s\}$ and $\{\lambda_s\}$ are both uniformly bounded in $V$ and $E$, respectively.}

As concerns the sequence $\{u_s\}$,
if it is unbounded, then, without any loss of
generality, we may assume that $\|u_s\|_V\to+\infty$, as $s\to \infty$. Inserting $v=0_V$ and $\rho=0_E$ into (\ref{eqnss12}), it has
\begin{eqnarray*}
&&\langle Au_s,u_s\rangle\le b(0_V,\lambda_s)-b(u_s,0_E) +J^0(\gamma u_s;-\gamma u_s)+ \langle f,u_s\rangle\\[2mm]
&&\qquad\le  \alpha_J
+\beta_J\|\gamma\|_{\mathcal L(V,X)}^\theta\|u_s\|_V^\theta
+\|f\|_{V^*}\|u_s\|_V,
\end{eqnarray*}
and
\begin{eqnarray*}
&&\frac{\langle Au_s,u_s\rangle}{\|u_s\|_V^{\max\{\theta,1\}}}\le \frac{\alpha_J}{\|u_s\|_V^{\max\{\theta,1\}}}+\frac{\beta_J\|\gamma\|_{\mathcal L(V,X)}^\theta}{\|u_s\|_V^{\max\{0,1-\theta\}}}+\frac{\|f\|_{V^*}}{\|u_s\|_V^{\max\{0,\theta-1\}}}.
\end{eqnarray*}
Now, the coercivity condition $H(A)$(iii) concludes a contradiction, so the sequence $\{u_s\}$ is uniformly bounded in $V$.

It remains to show that the sequence $\{\lambda_s\}$ is uniformly bounded
in $E$ too.
For any $w\in V\setminus\{0_V\}$, inserting $v= u_s-\frac{w}{\|w\|_V}$ and $\rho =\lambda_s$ into (\ref{eqnss12}), we have
\begin{eqnarray*}
&&b(\frac{w}{\|w\|_V},\lambda_s)\le
\langle Au_s,-\frac{w}{\|w\|_V}\rangle  +J^0(\gamma u_s;-\gamma \frac{w}{\|w\|_V} )+\langle f,\frac{w}{\|w\|_V}\rangle\\[2mm]
&&\quad=\langle Au_s,-\frac{w}{\|w\|_V} \rangle+\max_{\xi_{u_s}\in\partial J(\gamma u_s)}\langle\xi_{u_s},-\gamma \frac{w}{\|w\|_V} \rangle_{X^*\times X}+\langle f,\frac{w}{\|w\|_V} \rangle\\[2mm]
&&\qquad\le \|Au_s\|_{V^*}+\|\gamma^*\partial J(\gamma u_s)\|_{V^*}+\|f\|_{V^*}.
\end{eqnarray*}
Passing to supremum with
$w \in V\setminus\{0_V\}$, by using inequality (\ref{inequalityb}), we deduce
\begin{equation*}
\alpha_b\|\lambda_s\|_E\le \limsup_{w\in V\setminus\{0_V\}}b(\frac{w}{\|w\|_V},\lambda_s)\le \|Au_s\|_{V^*}+\|\gamma^*\partial J(\gamma u_s)\|_{V^*}+\|f\|_{V^*}
\end{equation*}
for all $s\in\mathbb N$.
This inequality combined with the uniform boundedness of $\{u_s\}$, hypotheses $H(A)$(iv) and $H(J)$(iii), implies that the sequence $\{\lambda_s\}$ is uniformly bounded in $E$.

\smallskip

From Claim 3, we can find $s_0\in\mathbb N$ large enough such that $\|\lambda_{s_0}\|_E<s_0$.
We shall verify that $(u_{s_0},\lambda_{s_0})$
is also a solution to Problem~\ref{problem1}.
Let $\eta\in \Lambda$ be arbitrary and $t\in(0,1)$ be small enough such that $\rho_t=t\eta+(1-t)\lambda_{s_0}\in Y(s_0)$. After inserting $\rho=\rho_t$ and $v=u_{s_0}$ into (\ref{eqnss12}) for $s=s_0$,
we have
\begin{equation*}
b(u_{s_0},\lambda_{s_0})-tb(u_{s_0},\eta)-(1-t)b(u_{s_0},\lambda_{s_0})\ge 0,
\end{equation*}
i.e.,
\begin{equation*}
b(u_{s_0},\eta-\lambda_{s_0})\le 0
\ \mbox{ for all }\ \eta\in \Lambda.
\end{equation*}
Putting $\rho=\lambda_{s_0}$ into (\ref{eqnss12}),
we have
\begin{equation*}
\langle Au_{s_0},v-u_{s_0}\rangle+b(v-u_{s_0},\lambda_{s_0}) +J^0(\gamma u_{s_0};\gamma v-\gamma u_{s_0})
\ge \langle f,v-u_{s_0}\rangle
\ \mbox{ for all } \ v\in V.
\end{equation*}
The last two inequalities reveal that $(u_{s_0},r_{s_0})\in V\times \Lambda$ is a solution to Problem~\ref{problem1},
which concludes the proof.
\hfill$\Box$

\medskip

In what follows, we denote the solution set to Problem~\ref{problem1} by $S(A,J,b,f)$. By Theorem~\ref{theorem2}, we know that the set $S(A,J,b,f) \subset V\times \Lambda$ is nonempty.
It is desirable to investigate further properties of the solution set which can be useful, for instance, in the study of optimal control problems for systems governed by a mixed variational-hemivariational inequality.
%
This is the reason that in the next theorem
we study essential properties of $S(A,J,b,f)$, such as, convexity, closedness and continuity.

\begin{Theorem}\label{theorem3}
Let $\Lambda$ be a nonempty, closed and convex subset of $E$ with $0_E\in \Lambda$. If hypotheses $H(A)$, $H(b)$, $H(J)$, $H(h)$(i)--(ii), and $H(\gamma)$ are satisfied, then the following hold

\smallskip

{\rm(i)} \
the solution set $S(A,J,b,f)$  is bounded and weakly closed in $V\times \Lambda$,

\smallskip

{\rm(ii)} \
if $h\colon V\to \real$ is convex, then the set $S(A,J,b,f)$ is convex as well,

\smallskip

{\rm(iii)} \
the multivalued mapping
$\mathcal S\colon V^*\to 2^{V\times\Lambda}$, defined by $f\mapsto \mathcal S(f):=S(A,J,b,f)$, is bounded, i.e., $\mathcal S$ maps bounded subsets of $V^*$ to bounded subsets
of $V\times \Lambda$,

\smallskip

{\rm(iv)} \ the multivalued mapping $f\mapsto \mathcal S(f)$ is strongly-weakly upper semicontinuous, (i.e., it is upper semicontinuous from $V^*$ endowed with the norm topology to the subsets of $V\times \Lambda$ endowed with the weak topology), and it has a strongly-weakly closed graph.
\end{Theorem}

\noindent
{\bf Proof.}
(i)
Arguing by contradiction, if we assume that $S(A,J,b,f)$ is unbounded, then there exists
a sequence $\{(u_n,\lambda_n)\}\subset S(A,J,b,f)$ such that
\begin{equation}\label{eqnss18}
\|u_n\|_V+\|\lambda_n\|_E\to +\infty,
\ \mbox{ as }\ n\to\infty.
\end{equation}
For the sequence $\{(u_n,\lambda_n)\}$, we claim that $\{u_n\}$ is bounded in $V$. If it is not true, then one has
\begin{equation}\label{eqnss19}
\|u_n\|_V \to +\infty,
\ \mbox{ as }\ n\to\infty.
\end{equation}
In fact, for each $n\in \mathbb N$, we have
\begin{eqnarray}
&&
\hspace{-1.2cm}
\langle A(u_n),v-u_n\rangle+b(v-u_n,\lambda_n)+J^0(\gamma u_n;\gamma v-\gamma u_n) \ge \langle f,v-u_n\rangle\mbox{ for all }v\in V,\label{eqnss20}\\[2mm]
&&
\hspace{-1.2cm}
b(u_n,\rho-\lambda_n)\le 0
\ \ \ \mbox{\rm for all}\ \
\rho\in \Lambda.\label{eqnss21}
\end{eqnarray}
Choosing $v=0_V$ and $\rho=0_E$ in (\ref{eqnss20}) and (\ref{eqnss21}), respectively, from the resulting inequalities,
it yields
\begin{eqnarray*}
&&\langle A(u_n),u_n\rangle\le b(-u_n,\lambda_n)+J^0(\gamma u_n;-\gamma u_n)+ \langle f,u_n\rangle\\[2mm]
&&\quad \le \alpha_J+\beta_J\|\gamma\|_{\mathcal L(V,X)}^\theta\|u_n\|_V^\theta+\|f\|_{V^*}\|u_n\|_V,
\end{eqnarray*}
and hence
\begin{equation*}
\frac{\langle A(u_n),u_n\rangle}{\|u_n\|_V^{\max\{1,\theta\}}}\le \frac{\alpha_J}{\|u_n\|_V^{\max\{1,\theta\}}}+\frac{\beta_J\|\gamma\|_{\mathcal L(V,X)}^\theta}{\|u_n\|_V^{\max\{1-\theta,0\}}}+\frac{\|f\|_{V^*}}{\|u_n\|_V^{\max\{0,\theta-1\}}}.
\end{equation*}
The above inequality combined
with (\ref{eqnss19}) and hypothesis
$H(A)$(iii) leads to a contradiction.
Hence, we deduce that the sequence $\{u_n\}$ is bounded in $V$.

Moreover, by (\ref{eqnss18}), we know that
$\{\lambda_n\}$ is unbounded in $E$.
Let $w\in V\setminus\{0_V\}$ be arbitrary.  Taking account of $v= u_n-\frac{w}{\|w\|_V}$
in (\ref{eqnss20}), it reads
\begin{eqnarray*}
&&b(\frac{w}{\|w\|_V},\lambda_n)\le
\langle Au_n,-\frac{w}{\|w\|_V}\rangle  +J^0(\gamma u_n;-\gamma \frac{w}{\|w\|_V} )+\langle f,\frac{w}{\|w\|_V}\rangle\nonumber\\[2mm]
&&\quad=\langle Au_n,-\frac{w}{\|w\|_V} \rangle+\max_{\xi_{u_n}\in\partial J(\gamma u_n)}\langle\xi_{u_n},-\gamma \frac{w}{\|w\|_V} \rangle_{X^*\times X}+\langle f,\frac{w}{\|w\|_V} \rangle\nonumber\\[2mm]
&&\qquad\le \|Au_n\|_{V^*}+\|\gamma^*\partial J(\gamma u_n)\|_{V^*}+\|f\|_{V^*}.
\end{eqnarray*}
Passing to supremum with
$w \in V\setminus\{0_V\}$ and using hypotheses $H(A)$(iv), $H(J)$(iii), $H(b)$, and the boundedness of $\{u_n\}$, we are able to find a constant $C^*$, which is independent of $n$, such that
\begin{equation*}
\|\lambda_n\|_E\le C^*
\ \mbox{ for all }\ n\in\mathbb N.
\end{equation*}
This leads to a contradiction with (\ref{eqnss18}). Therefore, we conclude that  the solution set $S(A,J,b,f)$ to Problem~\ref{problem1} is bounded
in $V\times \Lambda$.

Next, we show the weak closedness
of $S(A,J,b,f)$.
Let $\{(u_n,\lambda_n)\}\subset S(A,J,b,f)$ be a sequence such that
\begin{equation}
u_n\to u\mbox{ weakly in $V$, and $\lambda_n\to \lambda$ weakly in $E$, as $n\to \infty$}.
\end{equation}
It follows from Theorem~\ref{theorem1} that
\begin{eqnarray*}
&&
\hspace{-0.5cm}
\langle Au_n,v-u_n\rangle
+b(v,\lambda_n)-b(u_n,\rho)+J^0(\gamma u_n;\gamma v-\gamma u_n)\\[2mm]
&&\ge \langle f,v-u_n\rangle
\ \ \mbox{\rm for all}\
v\in K\ \mbox{\rm and all}\ \rho\in Y.
\end{eqnarray*}
Furthermore, the $h$-relaxed monotonicity of the mapping
$u\mapsto A u +\gamma^*\partial J(\gamma u)$ implies
\begin{eqnarray*}
&&
\hspace{-0.5cm}
\langle Au_n,v-u_n\rangle+J^0(\gamma u_n;\gamma v-\gamma u_n)+h(v-u_n)\\[2mm]
&&=\langle Au_n,v-u_n\rangle+\langle\xi_{u_n},\gamma (v-u_n)\rangle_{X^*\times X}+h(v-u_n)\\[2mm]
&&\le \langle Av,v-u_n\rangle+\langle\xi_v,\gamma (v-u_n)\rangle_{X^*\times X}
\end{eqnarray*}
for all $\xi_v\in \partial J(\gamma v)$, where $\xi_{u_n}\in \partial J(\gamma u_n)$ is such that
$\langle\xi_{u_n},\gamma(v-u_n)\rangle_{X^*\times X}= J^0(\gamma u_n;\gamma v-\gamma u_n)$.
We use the last two inequalities to obtain
\begin{eqnarray*}
&&
\hspace{-0.5cm}
\langle Av,v-u_n\rangle
+b(v,\lambda_n)-b(u_n,\rho)+\langle\xi_v,\gamma (v-u_n)\rangle_{X^*\times X}\\[2mm]
&&\ge \langle f,v-u_n\rangle+h(v-u_n)\ \ \mbox{for all }\xi_v\in \partial J(\gamma v),
\end{eqnarray*}
and all $(v,\rho)\in V\times \Lambda$.
Passing to the upper limit, as $n\to \infty$,
one has
\begin{eqnarray*}
&&
\hspace{-0.5cm}
\langle f,v-u\rangle+h(v-u)\\[2mm]
&&\le\limsup_{n\to\infty}\big[ \langle f,v-u_n\rangle+h(v-u_n)\big]\\[2mm]
&&\le \limsup_{n\to\infty}\big[\langle Av,v-u_n\rangle
+b(v,\lambda_n)-b(u_n,\rho)+\langle\xi_v,\gamma (v-u_n)\rangle_{X^*\times X}\big]\\[2mm]
&&\le \langle Av,v-u\rangle
+b(v,\lambda)-b(u,\rho)+\langle\xi_v,\gamma (v-u)\rangle_{X^*\times X}\ \mbox{ (for all $\xi_v\in \partial J(\gamma v)$)}\\[2mm]
&&\le \langle Av,v-u\rangle
+b(v,\lambda)-b(u,\rho)+J^0(\gamma v;\gamma v-\gamma u)
\end{eqnarray*}
for all $v\in V$ and $\rho\in \Lambda$.
We now invoke Theorem~\ref{theorem1} again to reveal that $(u,\lambda)\in S(A,J,b,f)$.
Hence, $S(A,J,b,f)$ is a weakly closed set.

\smallskip

\noindent
(ii) Assume that $h$ is a convex function.
Let $(u_1,\lambda_1)$,
$(u_2,\lambda_2)\in S(A,J,b,f)$ and $t\in(0,1)$. From Theorem~\ref{theorem1}, we have
\begin{eqnarray*}
&&
\hspace{-0.5cm}
\langle Au_i,v-u_i\rangle
+b(v,\lambda_i)-b(u_i,\rho)+J^0(\gamma u_i;\gamma v-\gamma u_i)\\[2mm]
&&\ge  \langle f,v-u_i\rangle
\end{eqnarray*}
for all $v\in V$ and $\rho\in \Lambda$.
Hence, for $i=1$, $2$, we get
\begin{eqnarray*}
&&
\hspace{-0.5cm}
\langle Av,v-u_i\rangle
+b(v,\lambda_i)-b(u_i,\rho)+\langle\xi_v,\gamma(v-u_i)\rangle_{X^*\times X}\\[2mm]
&&\ge  \langle f,v-u_i\rangle+h(v-u_i)
\end{eqnarray*}
for all $\xi_v\in \partial J(\gamma v)$ and  $(v,\rho)\in V\times \Lambda$.
Here, we have applied the $h$-relaxed monotonicity of
$A\cdot +\gamma^*\partial J(\gamma\cdot)$.
Denote $u_t=tu_1+(1-t)u_2$ and $\lambda_t=t\lambda_1+(1-t)\lambda_2$.
It follows that
\begin{eqnarray*}
&&
\hspace{-0.5cm}
\langle Av,v-u_t\rangle
+b(v,\lambda_t)-b(u_t,\rho)+\langle\xi_v,\gamma(v-u_t)\rangle_{X^*\times X}\\[2mm]
&&=t\big[\langle Av,v-u_1\rangle
+b(v,\lambda_1)-b(u_1,\rho)+\langle\xi_v,\gamma(v-u_1)\rangle_{X^*\times X}\big]\\[2mm]
&&\quad+(1-t)\big[\langle Av,v-u_2\rangle
+b(v,\lambda_2)-b(u_2,\rho)+\langle\xi_v,\gamma(v-u_2)\rangle_{X^*\times X}\big]\\[2mm]
&&\ge t\langle f,v-u_1\rangle+(1-t)\langle f,v-u_2\rangle+th(v-u_1)+(1-t)h(v-u_2)\\[2mm]
&&\ge  \langle f,v-u_t\rangle+h(v-u_t)
\end{eqnarray*}
for all $\xi_v\in \partial J(\gamma v)$ and  $(v,\rho)\in V\times \Lambda$, where the last inequality is obtained by using the convexity of function $h$.
Combining this inequality with Theorem~\ref{theorem1} and the fact
\begin{eqnarray*}
&&
\hspace{-0.5cm}
\langle Av,v-u_t\rangle
+b(v,\lambda_t)-b(u_t,\rho)+J^0(\gamma v;\gamma v-\gamma u_t)\\[2mm]
&&\ge\langle Av,v-u_t\rangle
+b(v,\lambda_t)-b(u_t,\rho)+\langle\xi_v,\gamma(v-u_t)\rangle_{X^*\times X}\ \mbox{ (for all $\xi_v\in \partial J(\gamma v)$)}\\[2mm]
&&\ge  \langle f,v-u_t\rangle+h(v-u_t)\ \,\mbox{ for all }v\in V\mbox{ and }\rho\in \Lambda,
\end{eqnarray*}
implies that $(u_t,\lambda_t)$
is also a solution to Problem~\ref{problem1}.
This proves that $S(A,J,b,f)$ is a convex set.

\smallskip

\noindent
(iii) If $\mathcal S\colon V^*\to 2^{V\times \Lambda}$ is not a bounded mapping,
then we are able to find a bounded set
$B\subset V^*$, and sequences
$\{f_n\}\subset B$,
$\{(u_n,\lambda_n)\}\subset V\times \Lambda$ with
$(u_n,\lambda_n)\in \mathcal S(f_n)$
for all $n\in \mathbb N$ such that $\|u_n\|_V+\|\lambda_n\|_E\to +\infty$, as $n\to\infty$.
As in the proof of assertion (i),
we derive
\begin{equation*}
\frac{\langle A(u_n),u_n\rangle}{\|u_n\|_V^{\max\{1,\theta\}}}\le \frac{\alpha_J}{\|u_n\|_V^{\max\{1,\theta\}}}+\frac{\beta_J\|\gamma\|_{\mathcal L(V,X)}^\theta}{\|u_n\|_V^{\max\{1-\theta,0\}}}+\frac{\|f_n\|_{V^*}}{\|u_n\|_V^{\max\{0,\theta-1\}}},
\end{equation*}
and
\begin{equation*}
b(\frac{w}{\|w\|_V},\lambda_n) \le \|Au_n\|_{V^*}+\|\gamma^*\partial J(\gamma u_n)\|_{V^*}+\|f_n\|_{V^*}.
\end{equation*}
The above inequalities combined with the boundedness of $\{f_n\}$ and the coercivity condition $H(A)$(iii) leads to a contradiction. Consequently, $\mathcal S\colon V^*\to 2^{V\times \Lambda}$ is a bounded mapping.

\smallskip

\noindent
(iv) In order to prove the upper semicontinuity of the mapping $f\mapsto \mathcal S(f)$, by  Proposition~\ref{posusc}, it is enough to prove that for each weakly closed subset $C$ of $V\times \Lambda$,
the set
\begin{equation*}
\mathcal S^-(C):=\big\{f\in V^* \mid \mathcal S(f)\cap C\neq\emptyset\big\}
\end{equation*}
is closed in $V^*$.
Let $C$ be a weakly closed subset
of $V\times \Lambda$, and
$\{f_n\}\subset \mathcal S^-(C)$ be a sequence such that $f_n\to f$ in $V^*$.
Thus, for each $n\in \mathbb N$,
we can find a pair of elements $(u_n,\lambda_n)\in \mathcal S(f_n)\cap C$, i.e.,
\begin{eqnarray*}
&&
\hspace{-0.5cm}
\langle Au_n,v-u_n\rangle
+b(v,\lambda_n)-b(u_n,\rho)+J^0(\gamma u_n;\gamma v-\gamma u_n)\\[2mm]
&&\ge  \langle f_n,v-u_n\rangle\ \mbox{ for all }v\in V\mbox{ and all }\rho\in \Lambda.
\end{eqnarray*}
The $h$-relaxed monotonicity of $A\cdot +\gamma^*\partial J(\gamma \cdot)$ shows that
\begin{eqnarray}\label{cv1}
&&
\hspace{-0.5cm}
\langle Av,v-u_n\rangle
+b(v,\lambda_n)-b(u_n,\rho)+\langle \xi_v,\gamma (v- u_n)\rangle_{X^*\times X}\nonumber\\[2mm]
&&\ge  \langle f_n,v-u_n\rangle+h(v-u_n)
\end{eqnarray}
for all $\xi_v\in \partial J(\gamma v)$ and  $(v,\rho)\in V\times \Lambda$.
From assertion (iii), we can see that the sequence $\{(u_n,\lambda_n)\}$
is bounded in $V\times E$.
The latter and the reflexivity of $V\times E$ guarantee that there exist a subsequence of $\{(u_n,\lambda_n)\}$,
still denoted in the same way,
and a pair of elements $(u,\lambda)\in V\times \Lambda$ such that
\begin{equation}\label{cv2}
u_n\to u\mbox{ weakly in $V$, and $\lambda_n\to \lambda$ weakly in $E$, as $n\to\infty$}.
\end{equation}
Taking into account the inequality (\ref{cv1}), and passing to the upper limit, as $n\to \infty$, we obtain
\begin{eqnarray*}
&&
\hspace{-0.5cm}
\langle Av,v-u\rangle
+b(v,\lambda)-b(u,\rho)+ J^0(\gamma v;\gamma v- \gamma u)\\[2mm]
&&\ge\langle Av,v-u\rangle
+b(v,\lambda)-b(u,\rho)+\langle \xi_v,\gamma (v- u)\rangle_{X^*\times X}\ \mbox{ (for all $\xi_v\in \partial J(\gamma v)$)}\\[2mm]
&&\ge  \langle f,v-u\rangle+h(v-u)
\end{eqnarray*}
for all $(v,\rho)\in V\times \Lambda$.
It follows from Theorem~\ref{theorem1} that $(u,\lambda)\in \mathcal S(f)$.
On the other hand, the convergences (\ref{cv2}) entail that $(u,\lambda)\in \mathcal S(f)\cap C$, due to the weak closedness of the set $C$. Therefore, we conclude that $\mathcal S$ is a strongly-weakly upper semicontinuous mapping.

Finally, since $\mathcal S$ is a strongly-weakly u.s.c. mapping with weakly closed values,
we are now in a position to apply Theorem~\ref{closedthem} to obtain the desired result that $\mathcal S$ has a strongly-weakly closed graph. This completes the proof.
\hfill$\Box$

\medskip

\begin{Remark}\label{remark2}
From the proofs of Theorems~\ref{theorem2} and~\ref{theorem3}, we can see that
the essence of the coercivity condition $H(A)$(iii) and inequality $H(J)$(ii) is to guarantee the following condition
\begin{equation}\label{coercivecondition}
\frac{\langle Au,u\rangle
	-J^0(\gamma u;-\gamma u)}{\|u\|_V}\to +\infty, \ \mbox{ as }\ \|u\|_V\to \infty.
\end{equation}
Moreover, it can be observed that if
$h$ is coercive in the following sense
\begin{equation}\label{hcoercive}
\frac{h(u)}{\|u\|_V}\to +\infty,
\ \mbox{ as }\ \|u\|_V\to \infty,
\end{equation}
and hypotheses $H(J)$(i), $H(J)$(iii), and $H(A)$(ii) hold, then condition (\ref{coercivecondition})
is automatically satisfied.
In that case, assumptions $H(J)$(ii) and $H(A)$(iii) could be removed.
\end{Remark}

\begin{Lemma}\label{lemma1}
Assume that $H(J)$(i), $H(J)$(iii), and $H(A)$(ii)  are fulfilled. If the function $h$ is coercive in the sense of (\ref{hcoercive}),
then condition (\ref{coercivecondition}) holds.
\end{Lemma}

\noindent
{\bf Proof.} Let $u\in V$.
By the $h$-relaxed monotonicity of $A\cdot+\gamma^*\partial J(\gamma\cdot)$,
we deduce
\begin{eqnarray*}
&&h(u)\le \langle Au-A0_V,u\rangle+\langle \xi_u-\xi_0,\gamma u \rangle_{X^*\times X}\ \mbox{(for all $\xi_0\in\partial J(0_X)$)}\\[2mm]
&&\quad\le \langle Au,u\rangle+\big(\|A0_V\|_{V^*}+\|\gamma^*\partial J( 0_X)\|_{V^*}\big)\|u\|_V-J^0(\gamma u;-\gamma u),
\end{eqnarray*}
where $\xi_u\in \partial J(\gamma u)$
is such that
$\langle\xi_u,-\gamma u
\rangle_{X^*\times X}
=J^0(\gamma u;-\gamma u)$.
This implies
\begin{equation*}
\frac{h(u)}{\|u\|_V}-\big(\|A0_V\|_{V^*}
+\|\gamma^*\partial J( 0_X)\|_{V^*}\big)
\le \frac{\langle Au,u\rangle
-J^0(\gamma u;-\gamma u)}{\|u\|_V}.
\end{equation*}
Now, by the coercivity condition (\ref{hcoercive}) and hypothesis $H(J)$(iii) we obtain the desired conclusion.
\hfill$\Box$

\begin{Remark}\label{remark3}
Note that if $A\colon V\to V^*$ is strongly monotone with constant $m_A>0$ and
$\partial J\colon X\to X^*$ is relaxed monotone with constant $m_J\ge 0$, i.e.,
\begin{equation*}
\langle Au-Av,u-v\rangle\ge
m_A\|u-v\|^2_V\ \ \ \mbox{\rm and}
\ \ \
\langle \xi_w-\xi_x,w-x \rangle_{X^*\times X}\ge -m_J\|w-x\|_X^2
\end{equation*}
for all $u$, $v\in V$, all $x$, $w\in X$,
and all $\xi_w\in\partial J(w)$,
$\xi_x\in \partial J(x)$, and
the inequality  $m_J\|\gamma\|_{\mathcal L(V,X)}^2<m_A$ holds,
then $h\colon V\to \real$ chosen as $h(u)=
(m_A-m_J\|\gamma\|_{\mathcal L(V,X)}^2)\|u\|_V^2$
satisfies condition (\ref{coercivecondition}).
\end{Remark}

In what follows, we introduce a multivalued function $\mathcal S_1\colon V^*\to 2^V$ defined by
\begin{eqnarray*}
\mathcal S_1(f):=\big\{u\in V\,\mid\,\mbox{ there exists }\lambda \in \Lambda \mbox{ such that }(u,\lambda)\in \mathcal S(f)\big\}.
\end{eqnarray*}

The unique solvability is of fundamental importance in numerical analysis
of the problem. So, this brings about the natural question of whether the mixed variational-hemivariational inequality has a unique solution. The following theorem examines a significant conclusion that Problem~\ref{problem1} has at least a solution $(u,\lambda)\in V\times \Lambda$, which is unique in its first component.
\begin{Theorem}\label{theorem4}
Let $\Lambda$ be a nonempty, closed and convex subset of $E$ with $0_E\in \Lambda$. If hypotheses $H(A)$, $H(b)$, $H(J)$, $H(h)$, and $H(\gamma)$ are fulfilled, then Problem~\ref{problem1} has at least one solution $(u,\lambda)\in V\times \Lambda$, which is unique in its first component.
\end{Theorem}

\noindent
{\bf Proof.}
The existence of solution is a direct consequence of Theorem~\ref{theorem2}.
Now, we shall prove the uniqueness in the first component of the pair solution for Problem~\ref{problem1}.
Let $(u_1,\lambda_1)\in V\times \Lambda$ and $(u_2,\lambda_2)\in V\times \Lambda$ be
solutions to Problem~\ref{problem1},
so, we have
\begin{eqnarray}
&&
\hspace{-1.2cm}
\langle A(u_i),v-u_i\rangle+b(v-u_i,\lambda_i)+J^0(\gamma u_i;\gamma v-\gamma u_i) \ge \langle f,v-u_i\rangle  \ \ \mbox{\rm for all}\ \ v\in V,\label{eqnss14} \\[2mm]
&&
\hspace{-1.2cm}
b(u_i,\rho-\lambda_i)\le 0
\ \ \ \mbox{\rm for all}\ \
\rho\in \Lambda \label{eqnss15}
\end{eqnarray}
with $i=1$, $2$.
We take $\rho = \lambda_2$ and $\rho = \lambda_1$ in (\ref{eqnss15}) for $i = 1$ and $2$, respectively, then, we sum up the resulting inequalities to obtain
\begin{equation}\label{eqnss16}
b(u_1-u_2,\lambda_2-\lambda_1) =b(u_1,\lambda_2-\lambda_1)+b(u_2,\lambda_1-\lambda_2)\le 0.
\end{equation}
On the other hand, inserting $v = u_2$ in (\ref{eqnss14}) for $i = 1$ and $v = u_1$ in (\ref{eqnss14}) for $i = 2$, accordingly, and combining the resulting inequalities with (\ref{eqnss16}),
we can find elements $\xi_{u_1}\in \partial J(\gamma u_1)$ and $\xi_{u_2}\in \partial J(\gamma u_2)$ such that
\begin{equation*}
\left\{\begin{array}{lll}
J^0(\gamma u_1;\gamma u_2-\gamma u_1)=\langle\xi_{u_1},\gamma u_2-\gamma u_1\rangle_{X^*\times X}, \\[2mm]
J^0(\gamma u_2;\gamma u_1-\gamma u_2)=\langle\xi_{u_2},\gamma u_1-\gamma u_2\rangle_{X^*\times X},
\end{array}\right.
\end{equation*}
and
\begin{eqnarray*}
&&
\hspace{-0.5cm}
\langle A(u_1)-A(u_2),u_1-u_2\rangle-\big(\langle\xi_{u_1},\gamma u_2-\gamma u_1\rangle_{X^*\times X}+\langle\xi_{u_2},\gamma u_1-\gamma u_2\rangle_{X^*\times X}\big)\\[2mm]
&&=\langle A(u_1)-A(u_2),u_1-u_2\rangle-\big(J^0(\gamma u_1;\gamma u_2-\gamma u_1)+J^0(\gamma u_2;\gamma u_1-\gamma u_2)\big)\\[2mm]
&&\le b(u_2-u_1,\lambda_1)+b(u_1-u_2,\lambda_2)= b(u_1-u_2,\lambda_2-\lambda_1)\le 0.
\end{eqnarray*}
The latter combined with the $h$-relaxed monotonicity of
$A\cdot +\gamma^*\partial J(\gamma \cdot)$
implies
\begin{equation*}
h(u_1-u_2)\le \langle A(u_1)+\gamma^*\xi_{u_1}-A(u_2)-\gamma^*\xi_{u_2},u_1-u_2\rangle\le 0.
\end{equation*}
Recalling that $h(v)>0$ for all
$v\in V\backslash\{0_V\}$,
we conclude that $u_1=u_2$.
\hfill$\Box$

\medskip
In what follows, when the first component of  the pair solution to Problem~\ref{problem1} is unique, we have the following stability result.
\begin{Theorem}\label{theorem5}
Let $\Lambda$ be a nonempty, closed and convex subset of $E$ with $0_E\in \Lambda$.
If hypotheses $H(A)$, $H(b)$, $H(J)$, $H(h)$,
and $H(\gamma)$ hold, then the mapping
$\mathcal S_1\colon V^*\to V $ is weakly continuous, i.e., $f_n\to f$ in $V^*$ implies $\mathcal S_1(f_n)\to \mathcal S_1(f)$ weakly in $V $.
Moreover, if there exist $\tau>1$ and $c_h>0$ such that $h(u)\ge c_h\|u\|^\tau_V$
for all $u\in V$, then we have
\begin{eqnarray}\label{eqnss22}
\|\mathcal S_1(f_1)-\mathcal S_1(f_2)\|_V\le c_h^\frac{1}{\tau-1} \,
\|f_1-f_2\|_{V^*}^\frac{1}{\tau-1}.
\end{eqnarray}
\end{Theorem}

\noindent
{\bf Proof.} Let $\{f_n\}\subset V^*$ be a sequence such that $f_n\to f$ in $V^*$.
Let $(u_n,\lambda_n)\in V\times \Lambda$
be a solution to Problem~\ref{problem1} corresponding to $f_n$. It follows from the assertion (iii) of Theorem~\ref{theorem3} that the sequence $\{(u_n,\lambda_n)\}$ is bounded in $V\times \Lambda$. The reflexivity of $V\times E$ ensures that there exist a subsequence of $\{(u_n,\lambda_n)\}$, still denoted by the same symbol, and a pair of elements
$(u,\lambda)\in V\times \Lambda$ such that
\begin{equation}\label{eqnss23}
u_n\to u\mbox{ weakly in $V$, and $\lambda_n\to \lambda $ weakly in $E$, as $n\to\infty$}.
\end{equation}
We now claim that $(u,\lambda)$ is also a solution of Problem~\ref{problem1} associated with $f$.
Indeed, for each $n\in\mathbb N$, we use Theorem~\ref{theorem1} to get
\begin{equation*}
\langle A(u_n),v-u_n\rangle+b(v,\lambda_n)-b(u_n,\rho)
+J^0(\gamma u_n;\gamma v-\gamma u_n)
\ge \langle f_n,v-u_n\rangle
\end{equation*}
for all $v\in V$ and $\rho\in \Lambda$.
From the $h$-relaxed monotonicity of
$A\cdot +\gamma^*\partial J(\gamma\cdot)$,
we have
\begin{equation}\label{eqnss25}
\langle A(v),v-u_n\rangle+b(v,\lambda_n)-b(u_n,\rho)+\langle\xi_v,\gamma (v-u_n)\rangle_{X^*\times X}\ge \langle f_n,v-u_n\rangle
\end{equation}
for all $\xi_v\in \partial J(\gamma v)$
and all $(v,\rho)\in V\times \Lambda$.
Passing to the upper limit in (\ref{eqnss25}),
as $n\to\infty$,
and applying Theorem~\ref{theorem1},
we conclude that $(u,\lambda)$ is a solution to Problem~\ref{problem1} with respect to $f$.
Note that the first component of the pair solution for Problem~\ref{problem1} is unique, this confesses that every subsequence of $\{u_n\}$ converges weakly to the same limit $u$,
so, we deduce that the whole sequence $\{u_n\}$ converges weakly to $u$ in $V$, thus is, $\mathcal S_1(f_n)\to \mathcal S(f)$ weakly in $V$, as $n\to\infty$.

Furthermore, we assume that there exist $\tau>1$ and $c_h>0$ such that $h(u)\ge c_h\|u\|^\tau_V$ for all $u\in V$. Let $(u_i,\lambda_i)$ be a solution to Problem~\ref{problem1} corresponding to $f_i$, for $i=1$, $2$, respectively. So, we have $u_i=\mathcal S_1(f_i)$ for $i=1$, $2$.
By an easy calculation, it turns out that
\begin{eqnarray*}
&&
\langle A(\mathcal S_1(f_1))-A(\mathcal S_1(f_2)),\mathcal S_1(f_1)-\mathcal S_1(f_2)\rangle-\big(J^0(\gamma \mathcal S_1(f_1);\gamma \mathcal S_1(f_2)-\gamma \mathcal S_1(f_1))\\[2mm]
&&+J^0(\gamma \mathcal S_1(f_2);\gamma \mathcal S_1(f_1)-\gamma \mathcal S_1(f_2))\big) \le \langle f_1-f_2,\mathcal S_1(f_1)-\mathcal S_1(f_2)\rangle.
\end{eqnarray*}
Next, $H(A)$(ii) implies
\begin{eqnarray*}
&&c_h\|\mathcal S_1(f_1)-\mathcal S_1(f_2)\|^\tau_V\le h(\mathcal S_1(f_1)-\mathcal S_1(f_2))\le \langle f_1-f_2,\mathcal S_1(f_1)-\mathcal S_1(f_2)\rangle\\[2mm]
&&\quad\le \|f_1-f_2\|_{V^*}\|\mathcal S_1(f_1)-\mathcal S_1(f_2)\|_V.
\end{eqnarray*}
Consequently, we can readily derive the inequality (\ref{eqnss22}).
This completes the proof.
\hfill$\Box$

\medskip

From Theorems~\ref{theorem2} and~\ref{theorem3}, we have the following result.
\begin{Corollary}\label{corollay1}
Let $\Lambda$ be a nonempty, closed and convex subset of $E$ with $0_E\in \Lambda$. If hypotheses $H(A)$, $H(b)$, $H(J)$, and $H(\gamma)$ are satisfied with $h(u)=0$ for all $u\in V$, then we have

\smallskip

{\rm(i)}
the solution set $S(A,J,b,f)$ to Problem~\ref{problem1} is nonempty, bounded, weakly closed, and convex in $V\times \Lambda$,

\smallskip

{\rm(ii)}
the multivalued mapping $\mathcal S\colon V^*\to 2^{V\times\Lambda}$ defined by
$f\mapsto \mathcal S(f):=S(A,J,b,f)$
is bounded, i.e.,
$\mathcal S$ maps bounded subsets of $V^*$ to bounded subsets of $V\times \Lambda$,

\smallskip

{\rm(iii)}
the multivalued mapping $f\mapsto \mathcal S(f)$ is strongly-weakly upper semicontinuous, (i.e., it is upper semicontinuous from $V^*$ endowed with the norm topology to the subsets of $V\times \Lambda$ endowed with the weak topology),
and it has a strongly-weakly closed graph.
\end{Corollary}

Finally, invoking Remark~\ref{remark2}, Lemma~\ref{lemma1}, Theorems~\ref{theorem3},~\ref{theorem4} and~\ref{theorem5}, we obtain the following conclusion.
\begin{Corollary}\label{corollay2}
Let $\Lambda$ be a nonempty, closed and convex subset of $E$ with $0_E\in \Lambda$. If hypotheses $H(A)$(i)--(ii), $H(A)$(iv), $H(b)$, $H(J)$(i), $H(h)$(ii), $H(J)$(iii) and $H(\gamma)$ are fulfilled with
$h(u)\ge c_h\|u\|^q$
for all $u\in V$ and some $c_h>0$, $q> 1$,
then Problem~\ref{problem1}
has at least one solution  $(u,\lambda)\in V\times \Lambda$, which is unique in its first component. Moreover, the following inequality holds
\begin{equation*}
\|\mathcal S_1(f_1)-\mathcal S_1(f_2)\|_V\le c_h^\frac{1}{q-1}
\, \|f_1-f_2\|_{V^*}^\frac{1}{q-1}
\ \ \mbox{\rm for all} \ \ f_1, f_2 \in V^*.
\end{equation*}
\end{Corollary}


\begin{thebibliography}{99}
	
\bibitem{han1}M. Barboteu, K. Bartosz, W. Han, T. Janiczko, Numerical analysis of a hyperbo\-lic hemivariational inequality arising in dynamic contact, {\it SIAM J. Numer. Anal.} {\bf 53} (2015), 527--550.

\bibitem{bartosz} K. Bartosz, M. Sofonea, The Rothe method for variational-hemivariational inequalities with applications to contact mechanics, {\it SIAM J. Math. Anal.} {\bf 48} (2016),  861--883.

\bibitem{carl} S. Carl, D. Motreanu, Extremal solutions of quasilinear parabolic inclusions with generalized Clarke's gradient, {\it J. Differential Equations} {\bf 191} (2003),  206--233.

\bibitem{cf} F.H. Clarke, {\it Optimization and Nonsmooth Analysis}, Wiley, Interscience, New York, 1983.

\bibitem{matei2018narwa} M.C. Cojocaru, A. Matei, Well-posedness for a class of frictional contact models via mixed variational formulations, {\it Nonlinear Anal.} {\bf 47} (2019), 127--141.

\bibitem{costearadulescuJOGO2012} N. Costea, V. R\u{a}dulescu, Inequality problems of quasi-hemivariational type involving set-valued operators and a nonlinear term, {\it J. Global Optim.} {\bf 52} (2012), 743--756.

\bibitem{DMP1}
Z. Denkowski, S. Mig\'orski, N.S. Papageorgiou,
{\it An Introduction to Non\-li\-near Analysis: Theory}, Kluwer Academic/Plenum Publishers, Boston, Dordrecht, London, New York, 2003.

\bibitem{DMP2}
Z. Denkowski, S. Mig\'orski, N.S. Papageorgiou,
{\it An Introduction to Nonlinear Analysis: Applications}, Kluwer Academic/Plenum Publishers, Boston, Dordrecht, London, New York, 2003.

\bibitem{fk} K. Fan, Some properties of convex sets related to fixed point theorems, {\it Math. Ann.} {\bf 266} (1984), 519--537.

\bibitem{hanreddy} W. Han, B. Reddy, On the finite element method for mixed variational in\-equa\-li\-ties arising in elastoplasticity, {\it SIAM J. Numer. Anal.} {\bf 32} (1995), 1778--1807.

\bibitem{han2} W. Han, M. Sofonea, M. Barboteu, Numerical analysis of elliptic hemivariational inequalities, {\it SIAM J. Numer. Anal.} {\bf 55} (2017), 640--663.

\bibitem{han3} W. Han, M. Sofonea, D. Danan, Numerical analysis of stationary variational-he\-mi\-variational inequalities,  {\it Numer. Math.} {\bf 139} (2018), 563--592.



\bibitem{Matei}
A. Matei, A variational approach via bipotentials for a class
of frictional contact problems,
{\it Acta Appl. Math.} {\bf 134} (2014), 45--59.

\bibitem{matei2014rwa} A. Matei, An existence result for a mixed variational problem arising from Contact Mechanics, {\it Nonlinear Anal.} {\bf 20} (2014), 74--81.

\bibitem{matei2014mma} A. Matei, An evolutionary mixed variational problem arising from frictional contact mechanics, {\it Math. Mech. Solids} {\bf 19} (2014), 225--241.

\bibitem{matei2015} A. Matei, Two abstract mixed variational problems and aplications in contact mechanics, {\it Nonlinear Anal.} {\bf 22} (2015), 592--603.

\bibitem{matei2018} A. Matei, A mixed hemivariational-variational problem and applications, {\it Comput. Math. Appl.} (2018), doi: 10.1016/j.camwa.2018.08.068.

\bibitem{matei2010} A. Matei, R. Ciurcea, Contact problems for nonlinearly elastic materials:  weak solvability involving dual Lagrange multipliers, {\it The ANZIAM Journal} {\bf 52} (2010), 160--178.

\bibitem{matei2018aa} A. Matei, S. Sitzmann, K. Willner, B.I. Wohlmuth,  A mixed variational formulation for a class of contact problems in viscoelasticity, {\it Appl. Anal.} {\bf 97} (2018), 1340--1356.

\bibitem{matei2016ima} A. Matei, M. Sofonea,
A mixed variational formulation for a piezoelectric frictional contact problem, {\it IMA J. Appl. Math.} {\bf 82} (2016), 334--354.

\bibitem{migorskiochal2009} S. Mig\'{o}rski, A. Ochal, Quasi-static hemivariational inequality via vanishing acceleration approach, {\it SIAM J. Math. Anal.} {\bf 41} (2009), 1415--1435.

\bibitem{smo1}
S. Mig\'{o}rski, A. Ochal, M. Sofonea,
{\it Nonlinear Inclusions and Hemivariational Inequalities. Models and Analysis of Contact Problems}, Advances in Mechanics and Mathematics {\bf 26}, Springer, New York, 2013.

\bibitem{migorskizeng1} S. Mig\'{o}rski, S.D. Zeng, Hyperbolic hemivariational inequalities controlled by evolution equations with application to adhesive contact model, {\it Nonlinear Anal.} {\bf 43} (2018), 121--143.

\bibitem{migorskizeng2} S. Mig\'{o}rski, S.D. Zeng, Penalty and regularization method for variational-hemivariational inequalities with application to frictional contact, {\it Z. Angew. Math. Mech.} {\bf 98} (2018), 1503--1520.

\bibitem{migorskizeng3}  S. Mig\'{o}rski, S.D. Zeng, A class of differential hemivariational inequalities in Banach spaces, {\it J. Global Optim.} {\bf 72} (2018), 761--779.

\bibitem{motreanu1}  D. Motreanu, P.D. Panagiotopoulos, {\it Minimax Theorems and Qualitative Properties of the Solutions of Hemivariational Inequalities}, Nonconvex Optim. Appl., vol. {\bf 29}, Kluwer, Dordrecht, 1998.

\bibitem{pana1983} P.D. Panagiotopoulos, Nonconvex energy functions, hemivariational inequalities and substationary principles, {\it Acta Mech.} {\bf 42} (1983), 160--183.

\bibitem{Pa1} P.D. Panagiotopoulos,
{\em Inequality Problems in Mechanics and Applications\/}, Birkh\"{a}user, Boston, 1985.

\bibitem{Pa2}
P.D. Panagiotopoulos, {\it Hemivariational Inequalities, Applications in Mechanics and Engineering}, Springer-Verlag,
Berlin, 1993.

\bibitem{SM2}
M. Sofonea, A. Matei,
{\it Vatiational Inequalities with Applications. A study of Antiplane Frictional Contact Problems, Advances in Mechanics and Mathematics}, vol. {\bf 18}, Springer, New York, 2009.

\bibitem{SM}
M. Sofonea, A. Matei,
{\it Mathematical Models in Contact Mechanics},
London Ma\-the\-matical Society Lecture Note Series {\bf 398}, Cambridge University Press, 2012.

\bibitem{sofoneamatei2015jogo}
M. Sofonea, A. Matei, History-dependent mixed variational problems in contact mechanics, {\it J. Global Optim.} {\bf 61} (2015), 591--614.

\bibitem{SMHistory}
M. Sofonea, S. Mig\'orski,
{\it Variational-Hemivariational Inequalities with Applications},
Chapman \& Hall/CRC, Monographs and Research Notes in Mathematics,
Boca Raton, 2017.

\bibitem{tanghuang2013JOGO} G.J. Tang, N.J. Huang, Existence theorems of the variational-hemivariational inequalities, {\it J. Global Optim.} {\bf 56} (2013), 605--622.

\bibitem{xiaofanghuangJOGO2015} Y.B. Xiao, X.M. Yang, N.J. Huang, Some equivalence results for well-posedness of hemivariational inequalities, {\it J. Global Optim.} {\bf 61} (2015),   789--802.

\bibitem{Zeidler} E. Zeidler,
{\it Nonlinear Functional Analysis and Applications} II A/B, Springer,
New York, 1990.

\end{thebibliography}
\end{document}